\documentclass[a4paper, 12pt]{article}

\usepackage{graphicx}
\usepackage{url}
\usepackage{amsbsy,amssymb}
\usepackage{amsmath}
\usepackage{abstract}
\usepackage{hyperref}
\usepackage{color}
\usepackage{amsmath,amsthm}
\usepackage{mathtools}
\usepackage{xfrac}
\usepackage{faktor}

\usepackage[all]{xy}

\newtheorem{Proposition}{Proposition}
\newtheorem{Corollary}{Corollary}
\newtheorem{Theorem}{Theorem}
\newtheorem*{Proof}{Proof}

\newtheorem{Definition}{Definition}
\newtheorem{Remark}{Remark}
\newtheorem{Example}{Example}
\newtheorem{Algorithm}{Algorithm}

\newcommand{\dcov}{d^{\nabla}}

\begin{document}
%%%%%%%%%%%%%%%%%%%%%%%%%%%%%%%%%%%%%%%%%%%%%%%%%%%%%%%%%%%%%%%%%%%%%%%%%%%%%%%%%%%%%%%%%%%%%%%%%%

%%%%%%%%%%%%%%%%%%%%%%%%%%%%%%%%%%%%%%%%%%%%%%%%%%%%%%%%%%%%%%%%%%%%%%%%%%%%%%%
%title
%%%%%%%%%%%%%%%%%%%%%%%%%%%%%%%%%%%%%%%%%%%%%%%%%%%%%%%%%%%%%%%%%%%%%%%%%%%%%%%
{\LARGE\centering{\bf{Inverting covariant exterior derivative}}}
%%%%%%%%%%%%%%%%%%%%%%%%%%%%%%%%%%%%%%%%%%%%%%%%%%%%%%%%%%%%%%%%%%%%%%%%%%%%%%%%%
%authors
%%%%%%%%%%%%%%%%%%%%%%%%%%%%%%%%%%%%%%%%%%%%%%%%%%%%%%%%%%%%%%%%%%%%%%%%%%%%%%%%%

\begin{center}
\sf{Rados\l aw Antoni Kycia$^{1,2,a}$, Josef \v{S}ilhan$^{1}$}
\end{center}

\medskip
 \small{
\centerline{$^{1}$Masaryk University}
\centerline{Department of Mathematics and Statistics}
\centerline{Kotl\'{a}\v{r}sk\'{a} 267/2, 611 37 Brno, The Czech Republic}
\centerline{\\}
\centerline{$^{2}$Cracow University of Technology}
\centerline{Department of Computer Science and Telecommunications}
\centerline{Warszawska 24, Krak\'ow, 31-155, Poland}
\centerline{\\}

\centerline{$^{a}${\tt
kycia.radoslaw@gmail.com} }
}

%%%%%%%%%%%%%%%%%%%%%%%%%%%%%%%%%%%%%%%%%%%%%%%%%%%%%%%%%%%%%%%%%%%%%%%%%%%%%%
%absract
%%%%%%%%%%%%%%%%%%%%%%%%%%%%%%%%%%%%%%%%%%%%%%%%%%%%%%%%%%%%%%%%%%%%%%%%%%%%%%
\begin{abstract}
\noindent

The algorithm for inverting covariant exterior derivative is provided. It works for a sufficiently small star-shaped region of a fibered set - a local subset of a vector bundle and associated vector bundle. The algorithm contains some constraints that can fail, giving no solution, which is the expected case for parallel transport equations. These constraints are straightforward to obtain in the proposed approach. The relation to operational calculus and operator theory is outlined. The upshot of this paper is to show, using the linear homotopy operator of the Poincare lemma, that we can solve the covariant constant and related equations in a geometric and algorithmic way. The considerations related to the regularity of the solutions are provided.

\end{abstract}
Keywords: covariant exterior derivative; Poincar\'{e} lemma; antiexact forms; homotopy operator; fibered set\\
%MSC
Mathematical Subject Classification: 53-08, 53B05, 53Z05, 53B50 \\

%%%%%%%%%%%%%%%%%%%%%%%%%%%%%%%%%%%%%%%%%%%%%%%%%%%%%%%%%%%%%%%%%%%%%%%%%%%%
\section{Introduction}
%%%%%%%%%%%%%%%%%%%%%%%%%%%%%%%%%%%%%%%%%%%%%%%%%%%%%%%%%%%%%%%%%%%%%%%%%%%%

The covariant exterior derivative on associated vector bundles \cite{NaturalOperations} is one of the primary tools of modern differential geometry \cite{KobayashiNomizu, NaturalOperations}. Moreover, the fundamental physical theories as gauge theories (e.g., Electrodynamics, Strong and Weak interactions theories) are formulated within this framework \cite{Thirring, Bleecker, BennTucker, CovariantDerivativeInPhysics, Sternber, Rund, GagueTheories, ConnectionsInPhysics}. Therefore, a practical way of inverting (at least locally) the covariant exterior derivative is in demand.

We consider a fibered set $U\times V\rightarrow U$, where $U\in \mathbb{R}^{n}$ is a star-shaped subset, and $V$ is a vector space. Denote the space of $V$-valued differential forms on $U$ using $\Lambda(U,V)=\Lambda^{0}(U,V)\oplus \ldots \oplus\Lambda^{n}(U,V)$. The forms are sections of exterior power of the cotangent bundle of the fibered set. Then we can define the differential operator
\begin{equation}
 \dcov := d + A\wedge\_, \quad \dcov: \Lambda^{*}(U, V) \rightarrow  \Lambda^{*+1}(U, V)
\end{equation}
where $A\in \Lambda^{1}(U, End(V))$. This setup naturally extends to the local trivialization of a vector bundle $\pi: W\rightarrow M$, where we provide $\pi(U)=\tilde{U} \subset M$, and where $U$ is again assumed to be star-shaped for the possibility of application of the Poincare lemma.

For associated vector bundle $E = P\times_{G}V$ with the principal bundle $P\rightarrow M$, the Lie group $G$, and the vector space $V$, $\dcov$ is the exterior covariant derivative that acts on $\Lambda_{b}(P,V;\cdot)$ - basic (horizontal and equivariant) differential forms with the action of $G$ on $V$ is denoted by $\cdot$, and $A\in \Lambda^{1}(P, End(V); Ad)$ is an $Ad$-equivariant connection one-form \cite{KobayashiNomizu, NaturalOperations, Husemoller, Sternber}. By analogy, we will call $\dcov$ the covariant exterior derivative operator in any case considered in this paper. The structure of $\dcov$ is the obvious choice for the operator that raises the degree of a form by one. Moreover, $A$ will be called the connection form for the same reason. Similarly, as for the connection form for the associated vector bundle, we assume that $A$ has a constant rank.

The main result of the paper is an operational way of solving for $\phi\in \Lambda^{k}(U,V)$ the equation
\begin{equation}
 \dcov \phi = J,
\end{equation}
for a fixed form $J\in \Lambda^{k+1}(U,V)$, $0 \leq k < n$. These results will be presented in the central Section \ref{Section_Inversionformula} of the paper. The solution method involves constraints which, when it fails, leads to the conclusion that the solution does not exist. However, we can obtain the solution algorithmically, and the structure of the algorithm enables us to find the constraints effectively. We also provide a practical way of solving the curvature equation treated as the square of $\dcov$, i.e.,
\begin{equation}
 \left(\dcov\circ\dcov\right) \phi=J,
\end{equation}
where $J\in \Lambda^{k+2}(U,V)$ and $0 \leq k <n-1$. It is an algebraic equation. However, we propose a solution using the fact that it is the square of a differential operator. Moreover, by splitting higher order differential operators into the composition of covariant derivative and its Hodge dual, we can iteratively apply our method to solve more elaborated differential equations called here geometrically-based differential equations. We also want to stress-out that the solution $\phi$ will be only continuous with $d\phi$ continuous event for smooth $A$. The regularity issues will be also examined.

Our approach is not the first attempt to geometrically solve the homogenous parallel transport equation. To our knowledge, one existing practical and algorithmic way of inverting covariant exterior derivative on vector bundles is Chen's method of iterated integrals \cite{Chen, Chen2, Gugenheim-ChenIntergal}. The recent reformulation in the supermanifolds setup is presented in \cite{SuperPoincareLemma, SuperPoincareLemma2}. This idea is also used in perturbative Quantum Field Theory \cite{Kreimer-ChenIntergal}, or Yang-Mills theory \cite{PoincareLemmaForConnection}. The general idea behind Chen's iterated integrals method is to start from constructing Path Space for a given (super)vector bundle $\pi:V\rightarrow M$ with the connection one-form $\omega \in \Lambda^{1}(M,End(V))$. Then the parallel transport operator \cite{SuperPoincareLemma, SuperPoincareLemma2}, $\Phi(t):V_{\gamma(t=0)}\rightarrow V_{\gamma(t)}$, for a path $\gamma:[0;1]\rightarrow M$ is the solution of the following ODE:
\begin{equation}
 \frac{d\Phi^{\omega}(t)}{dt}= i_{t}^{*} \partial_{t} \lrcorner \omega\wedge \Phi^{\omega}(t), \quad \Phi^{\omega}(0)=Id_{V},
\end{equation}
where $i_{t}:M\rightarrow M\times[0,1]$ is the inclusion $i_{t}(x)=(x,t)$,  $\omega\in \Lambda^{*}(M\times [0,1], End(V))$ is a connection one-form on path space and $t \rightarrow \Phi^{\omega}(t)$ is a smooth mapping from the path parameter $t\in[0,1]$ to $\Lambda^{*}(M, End(V))$. Then the solution is the operator series (see, e.g., section 3 of \cite{SuperPoincareLemma})
\begin{equation}
 \Phi^{\omega}(t)=\sum_{n=0}^{\infty}\Phi^{\omega}_{n}(t),\quad \Phi^{\omega}(0)=Id_{V},
\end{equation}
where
\begin{equation}
 \begin{array}{c}
  \Phi^{\omega}_{0}(t)=Id, \\
  \Phi^{\omega}_{1}(t)=\int_{0}^{t}i^{*}_{s_{1}} \partial_{s_{1}} \lrcorner\omega ds_{1} \\
  \Phi^{\omega}_{n}(t)=\int_{0}^{t}\int_{0}^{s_{1}}\ldots\int_{0}^{s_{n-1}} i^{*}_{s_{1}} \partial_{s_{1}} \lrcorner\omega \wedge \ldots  \wedge i^{*}_{s_{n}} \partial_{s_{n}} \lrcorner\omega ~ ds_{n}\ldots ds_{1}, \quad n \geq 2,
 \end{array}
\end{equation}
for $t \geq s_{1} \geq \ldots \geq s_{n} \geq 0$. The formula relies on a path $\gamma$ and explores some basic ideas of the Poincare lemma. Chen noted that his integral operation is nilpotent (see the Corollary after Lemma 1.5.1 in \cite{Chen2}). We provide a method that explores this idea in depth using the calculus formulated by Edelen \cite{EdelenExteriorCalculus, EdelenIsovectorMethods, KyciaPoincare}.
In contrast with Chen's approach, our formula is valid in a local star-shaped subset of a bundle and not only along a path. Moreover, it can be applied to a more general setup of associated vector bundles and inhomogeneous equations. Chen's approach simplifies the PDE problem to the ODE problem defined along a path, and therefore additional construction of Path Space is needed. Our method solves general parallel transport PDE and therefore is more general. It can be classified as one of the methods of the vast discipline of Exterior Differential Systems \cite{Bryant, ExteriorDiffSys2}.

Our method relies on a simple observation that the covariant exterior derivative consists of the exterior derivative, $d$, and an algebraic operation that is the wedge product of a connection $1$-form. Then the exterior derivative can be locally inverted utilizing a homotopy operator of the Poincare lemma, and the connection form term introduces an additional level of complexity. The general approach to the homotopy operator is a classical subject; see, e.g., \cite{deRham}. However, the profits of using the homotopy operator defined for the linear homotopy \cite{EdelenExteriorCalculus, EdelenIsovectorMethods} was not widely recognized, although it proved to have many valuable properties that can be used to solve local problems in mathematics and mathematical physics, see \cite{KyciaPoincare, KyciaPoincareCohomotopy, CopoincareHamiltonianSystem} for examples and references.

Since we will use modern ideas related to the Poincare lemma, we will summarize the ideas behind this classical lemma for smooth differential forms, pointing out modern results used in our approach. The Poincar\'{e} lemma can be extended in many various directions \cite{deRham}, including non-abelian case \cite{NonabelianPoincareLemma}, higher order version \cite{HigherPoincare}, or be used in quantization \cite{PoincareQuantisation}, to mention a few. The Poincare lemma states, in the most practical formulation, that on a star-shaped subset $U \subset \mathbb{R}^{n}$ every closed form (an element of the kernel of exterior derivative $d$) is also exact (an element of the image of $d$). It can be proved in many ways, however, the usage of the linear homotopy operator is the most profitable one in our approach. To this end, let $x_{0}\in U$ be a center of the linear homotopy $F:U\times [0;1] \rightarrow U$, $F(x,t)=x_{0} + t(x-x_{0})$. Then the homotopy operator defined by $F$ is
\begin{equation}
H\omega = \int_{0}^{1} \mathcal{K}\lrcorner\omega|_{F(t,x)}t^{k-1}dt,
\end{equation}
for $\mathcal{K}=(x-x_{0})^{i}\partial_{i}$, and where $\omega \in \Lambda^{k}(U)$ is a $k$-differential form. It was noticed in \cite{EdelenExteriorCalculus, EdelenIsovectorMethods} that it is nilpotent\footnote{Nilpotency was also noticed by Chen, see Corollary after Lemma 1.5.1 in \cite{Chen2}, however not used to construct the inversion for covariant derivative.} $H^{2}=0$, and possesses useful properties like: $HdH=H$, $dHd=d$. However, the most useful is the homotopy invariance formula \cite{deRham, EdelenExteriorCalculus, EdelenIsovectorMethods, KyciaPoincare, KyciaPoincareCohomotopy}
\begin{equation}
 dH + Hd = I - s_{x_{0}}^{*},
 \label{Eq_homotopyFormula}
\end{equation}
where $s_{x_{0}}^{*}$ is the pullback along the constant map $s_{x_{0}}:x_{0} \hookrightarrow U$. The map $s_{x_{0}}^{*}$ can be nonzero only on $\Lambda^{0}(U)$.

Because of $(dH)^{2}=dH$ and $(Hd)^{2}=dH$, these operators can be interpreted as projection operators, and the formula (\ref{Eq_homotopyFormula}) as a decomposition of the identity operator to projectors into two spaces \cite{EdelenExteriorCalculus, EdelenIsovectorMethods}: $\Lambda^{*}(U)=\mathcal{E}(U)\oplus \mathcal{A}(U)$, where we denote an exact vector space (equivalent to closed forms on a star-shaped subset by the Poincare lemma) by $\mathcal{E}(U)=\{\omega \in \Lambda^{*}(U) | d\omega =0\}$, and the module over $C^{\infty}(U)$ of antiexact forms defined by $\mathcal{A}(U)=\{  \omega \in \Lambda^{*}(U) | \mathcal{K}\lrcorner \omega = 0, \omega|_{x=x_{0}}=0 \}$. It can be also shown that $\mathcal{E} = im(dH)$ and $\mathcal{A}=im(Hd)$. The direct sum decomposition is called the geometric decomposition.

Likewise, \cite{KyciaPoincareCohomotopy, CopoincareHamiltonianSystem}, for a star-shaped subset $U$ a Riemannian manifold $(M, g)$ with a metric tensor $g:TM\times TM \rightarrow \mathbb{R}$, we have the Hodge star $\star:\Lambda^{k}(M)\rightarrow \Lambda^{n-k}(M)$, $0\leq k \leq n=dim(M)$, that provides the codifferential $\delta= \star^{-1}d\star \eta$, for $\eta\omega = (-1)^{k}\omega$, where $\omega\in \Lambda^{k}$. In this setup, the dual theory of the co-Poincare lemma, and cohomotopy operator $h=\eta \star^{-1}H\star$ is defined. There exist (co)homotopy invariance formula \cite{KyciaPoincareCohomotopy, CopoincareHamiltonianSystem}, $h\delta+\delta h = I-S_{x_{0}}$, where $S_{x_{0}}=\star^{-1}s^{*}_{x_{0}}\star$. The cohomotopy operator $h$ defines the other geometric decomposition $\Lambda^{*}(U)=\mathcal{C}(U)\oplus \mathcal{Y}(U)$ into coexact (also coclosed) vector space $\mathcal{C}(U)=ker(\delta)=im(\delta h)$ and anticoexact module $\mathcal{Y}(U)=\{\omega \in \Lambda^{*}(U)| \mathcal{K}^{\flat}\wedge \omega=0\quad \omega|_{x=x_{0}}=0 \}=im(h\delta)$ over $C^{\infty}(U)$. Here $\flat:TM\rightarrow T^{*}M$ is a musical isomorphism related to the existence of a metric structure on $M$. The opposite isomorphism is given by $\sharp: T^{*}M\rightarrow TM$. These statements are easily obtained  from the (anti)exact theory by using the identity(e.g., \cite{BennTucker})
\begin{equation}
\alpha^{\sharp}\lrcorner \star\phi = \star(\phi\wedge\alpha),
\label{Eq.HodgeDuality}
\end{equation}
that dualizes (anti)exact theory to (anti)coexact one, see \cite{KyciaPoincareCohomotopy, CopoincareHamiltonianSystem}.

Both these decompositions into (co)(anti)exact forms of $\Lambda^{*}$ allow us to solve plenty of practical problems of mathematical physics, see \cite{KyciaPoincareCohomotopy, CopoincareHamiltonianSystem, KyciaPoincare} for modern applications, and \cite{EdelenExteriorCalculus, EdelenIsovectorMethods} for the classical ones. These results by a componentwise application can be extended to the vector-valued differential forms. In this case  we will denote corresponding spaces by $\mathcal{E}(U,V)$, $\mathcal{A}(U,V)$, $\mathcal{C}(U,V)$ and $\mathcal{Y}(U,V)$.

In this paper we use (co)(anti)exact decomposition to solve equations involving covariant exterior derivative in a star-shaped subset of a fibered set, that is, local trivialization of a general fiber bundle. We obtain a practical algorithm for getting covariantly constant differential forms. The paper is organized as follows: In the next section we provide the formula for inverting covariant exterior derivative. We also discuss the problem of applying our results to associated vector bundles. Next, we cast these formulas in Bittner's operational calculus framework, summarized for the reader's convenience in the Appendices \ref{Appendix_BittnersOperatorCalculus} and \ref{Appendix_BittnersOperatorCalculusAsCategory}. Then we use these results to provide an algorithm for solving the curvature equation. Next, we Hodge-dualize the previous results for manifolds with metrics structure. We also apply the geometric decomposition to the Cartan structure equations. Finally, we propose an Algorithm for the iterative solution of equations resulting from composing covariant derivatives and their Hodge dual operators. In the Appendices we discuss Bittner's operator calculus, the application of the results to the operator-valued connection, and the integral equations version of our results. We also motivate the introduction of the connection one-form when solving $d\phi=J$ in a similar, yet 'Cartan-like', manner as a 'minimal coupling' is used in physics for introducing a covariant derivative. In the final Appendix, some considerations related to the regularity of the solutions are given.

%%%%%%%%%%%%%%%%%%%%%%%%%%%%%%%%%%%%%%%%%%%%%%%%%%%%%%%%%%%%%%%%%%%%%%%%%%%%
\section{Inversion formula for covariant exterior derivative}
\label{Section_Inversionformula}
%%%%%%%%%%%%%%%%%%%%%%%%%%%%%%%%%%%%%%%%%%%%%%%%%%%%%%%%%%%%%%%%%%%%%%%%%%%%

The following two theorems provide the local inverse of the $\dcov$ operator. We want to underline that all considerations will be local, and the underlying set $U$ is an Euclidean star-shaped ones.

We must control the class of differentiability of (anti)(co) closed forms, and so we define:
\begin{itemize}
 \item {$\Lambda^{k,p}(U, V)$ - $k$-differenital froms with $C^{p}$ coefficients;}
 \item {$\mathcal{E}^{k,p}(U,V) = ker(d)\cap C^{p}(U,V)$,}
 \item {$\mathcal{A}^{k,p}(U,V) = ker(H)\cap C^{p}(U,V)$, }
 \item {$\mathcal{C}^{k,p}(U,V) = ker(\delta)\cap C^{p}(U,V)$, }
 \item {$\mathcal{Y}^{k,p}(U,V) = ker(h)\cap C^{p}(U,V)$.}
\end{itemize}
When index $p$ is omitted, it means that it is an union of all regularity classes for $p=\overline{0,\infty}$.

For simplicity we assume that $A\in \Lambda^{1,\infty}(U,End(V))$, however, for most considerations it suffices to have $A \in \Lambda^{1,1}(U,End(V))$ since only continuity of $A$ and $dA$ coefficients is used.

For convergence issues we must redefine all the notions to have more control over regularity of the coefficients of differential forms. We use open star-shaped subset $U \subset M$, that is bounded. Its closure is a compact set. Therefore we can define a norm suitable for uniform convergence considerations.
\begin{Definition}
Select the base $\{e_{i}\}_{i=1}^{dim(V)}$ of $V$. For a star-shaped subset of $ U$ with the center $x_{0}$.
For a $k$-form $\omega \in \Lambda^{*}(U,V)$ we have representation
$\omega = \omega^{i}_{I}dx^{I}e_{i}$, for a multiindex $I = (i_{1}, \ldots, i_{n})$, $i_{k}\in\{0,1\}$, where $dx^{I}=dx_{1}^{i_{1}}\wedge\ldots\wedge dx_{n}^{i_{n}}$. Then we define the norm
\begin{equation}
 ||\omega||_{\infty} = max_{i,I}sup_{x\in K} |\omega^{I}_{i}(x)|,
 \label{Eq.Norm_Vect_form}
\end{equation}
which is a simple extension of supremum norm for a vector-valued differential forms. Since $U$ is bounded so supremum is the maximum on $\bar{U}$.

In addition for an operator $A\wedge\_ = [A^{i}_{j,k}(x)dx^{k}]\wedge\_$ the supremum norm is given by
\begin{equation}
 ||A||_{\infty} = max_{ijk} sup_{x\in K} |A^{i}_{j,k}(x)|.
 \label{Eq.Norm_End_form}
\end{equation}
\end{Definition}
The above norms are extension of a norm for scalar-valued forms presented in \cite{deRham}. For examining higher regularity of the coefficients of the differential forms we must introduce supremum norm that involves derivatives. It will be done in Appendix \ref{Section_RegularityConsiderations}.

We will develop a complete theory starting from the homogenous parallel transport equation, then adding inhomogeneity as an exact form, and finally, pass to the general form of inhomogeneity.

%%%%%%%%%%%%%%%
\subsection{Homogenous equation}
%%%%%%%%%%%%%%%
First, we solve the homogenous equation of covariant constancy. We start from the scalar-valued differential forms. We have trivial:

%%%%%%%%%%%%%%%%%%%%%%%%%%%%%%%%%%
\begin{Proposition}
The unique smooth solution $\phi\in \Lambda^{0,\infty}(U,\mathbb{R})$ to the equation
 \begin{equation}
  \dcov \phi =0, \quad \phi\in \Lambda^{0, \infty}(U),~ A\in\Lambda^{1,\infty}(U,\mathbb{R}),
  \label{Eq_homogenous_equation}
 \end{equation}
 is given by
 \begin{equation}
  \phi = c\exp(-H(A)),
  \label{eq_solutionLambda0}
 \end{equation}
where $c \in \mathbb{R}$ is treated as an exact form (constant function) on $U$.
\end{Proposition}

\begin{Proof}
 We have
 \begin{equation}
  d\phi =-A\phi.
 \end{equation}
 We can easily check that the solution is $\phi = c\exp(-HA)$. For $c=0$, we obtain $\phi=0$, which is also a solution.
\end{Proof}

\begin{Remark}
Taking exterior derivative of (\ref{Eq_homogenous_equation}), we note that the solution $\phi$ is also a multiplier for $A$, i.e., $d(A\wedge \phi)=0$.

Moreover, note that the $c\in \mathbb{R}$ can be treated as a pullback of the solution (\ref{eq_solutionLambda0}) to the boundary point $x_{0}$, i.e., to the condition $\phi(x_{0})=c$.
\end{Remark}

Passing to vector-valued differential forms $\Lambda(U,V)$, we must first consider the $ker(A\wedge\_)$. We must deal with two cases since $\Lambda(U,V)\setminus ker(A\wedge\_)=\mathcal{E}(U,V)\setminus ker(A\wedge\_)\oplus \mathcal{A}(U,V) \setminus ker(A\wedge\_)$.

%%%%%%%%%%%%%%%%%%%%%%%%%%%%%%%%%%
\begin{Definition}
The solutions of $\dcov \phi=0$ fulfilling $\phi \in \mathcal{E}(U,V)\cap ker(A\wedge\_)$ we call the \textit{exact gauge modes}.
\end{Definition}

As we will see, the gauge modes are responsible for the nonuniqueness of the solution.

\begin{Proposition}
\label{Proposition_noantiexactGaugeModes}
 There is no solution to $\dcov \phi=0$ apart from the trivial one if we require that $\phi \in \mathcal{A}(U,V) \cap ker(A\wedge\_)$. We will call elements from $\mathcal{A}(U,V) \cap ker(A\wedge\_)$ \textit{antiexact gauge modes} if exist.
\end{Proposition}
\begin{Proof}
It there would be $\phi \in \mathcal{A}(U,V) \cap ker(A\wedge\_)$ then upon substituting into the equation we would have that $d\phi=0$, that is $\phi \in \mathcal{E}(U,V)\cap ker(A\wedge\_)$. The only element belonging to $\mathcal{E}\cap\mathcal{A}$ is $\phi=0$.
\end{Proof}

We also have obvious
\begin{Corollary}
If $\phi$ is a solution of $\dcov \phi=0$, then it is also a solution of $F\wedge \phi =0$, where $F =\dcov \dcov=dA + A\wedge A$ is the curvature of $A$
\end{Corollary}

We can now provide the solution to the full equation.
%%%%%%%%%%%%%%%%%%%%%%%%%%%%%%%%%%
\begin{Theorem}
\label{Th_homogenous_solution}
The unique nontrivial solution to the equation
 \begin{equation}
  \dcov \phi =0, ~ k>0
  \label{Eq_homogenous_equation}
 \end{equation}
 with the condition $dH\phi= c\in \mathcal{E}(U,V)\setminus ker(A\wedge\_)$, $c\neq 0$, is given by
 \begin{equation}
  \phi = \sum_{l=0}^{\infty} (-1)^{l} (H(A\wedge \_))^{l} c,
  \label{Eq.Solution_homogenous_k_gt_0}
 \end{equation}
where $c$ is an arbitrary form, $(H(A\wedge \_))^{0}=Id$, and
\begin{equation}
 (H(A\wedge \_))^{l} = \underbrace{H(A\wedge ( \ldots  (H(A\wedge \_ )\ldots )}_{l},
\end{equation}
is the $l$-fold composition of the operator $H\circ A \wedge \_$.

The series in (\ref{Eq.Solution_homogenous_k_gt_0}) is uniformly convergent to the continuous form ($\Lambda^{*,1}(U,V)$) inside the ball $B(x_{0}, r)$ with center $x_{0}$ and the radius $r$ given by
\begin{equation}
 r \frac{||A||_{\infty}}{k} = 1,
 \label{Eq_convergence_homogenous_solution}
\end{equation}
 where supremum is taken over the ball. Moreover, $d\phi = A\wedge \phi$ is a continuous form.

 The continuous solution is an element of $\Lambda^{k,0}(U,V)\setminus ker(A\wedge\_)$. Two of the solutions differ by an exact gauge mode obtained from (\ref{Eq.Solution_homogenous_k_gt_0}) when $c\in \mathcal{E}(U,V) \cap ker(A\wedge\_)$, i.e., $\phi=c$.
\end{Theorem}

The proof will use the perturbation series approach with decomposition into (anti)exact forms. This method can also be seen as a solution of the integral equation version of (\ref{Eq_homogenous_equation}) obtained by applying the $H$ operator and using (anti)exact decomposition -  see Appendix \ref{Appendix_Integral_equations}. Therefore some parts of the proof remind the derivation of the Neumann series for integral equations \cite{IntegralEquations}. It is analogous to the Picard-Lindel\"{o}f theorem for ODE.

%%%%%%%%%%%%%%%%%%%%%%%%%%%%%%%%%%
\begin{Proof}

We notice first, by taking exterior derivative of (\ref{Eq_homogenous_equation}), that $d(A\wedge \phi)=0$, i.e., $A\wedge \phi \in \mathcal{E}(U,V)$, i.e., $A\wedge \phi = dH(A\wedge\phi)$.

We split the proof into two main parts. The first one is formal derivation of the series (\ref{Eq.Solution_homogenous_k_gt_0}). The second part is devoted in proving that the formal series is uniformly convergent.

%%%%%
\textbf{Formal solution - differential equations approach:}
%%%%%

In order to introduce formal perturbation series, we modify the equation (\ref{Eq_homogenous_equation}) to
\begin{equation}
 d\phi + \lambda A\wedge \phi=0,
\end{equation}
introducing a real parameter $\lambda \neq 0$.

Searching the solution in the form of a formal power series
\begin{equation}
 \phi = \phi_{0} + \lambda \phi_{1} + \lambda^{2} \phi_{2} + \ldots,
\end{equation}
we get a system of equations with respect to the powers of $\lambda$. There will be some additional constraints for $\phi_{i}$ functions. We will assume that they are formally fulfilled at this point since we want to only heuristically motivate (\ref{Eq.Solution_homogenous_k_gt_0}) instead stating it without any comment. When we pass formally to the limit $\lambda=1$ these constraints will become irrelevant, since (\ref{Eq.Solution_homogenous_k_gt_0}) is much more general. It remains only the trivial constraint $F\wedge \phi=0$ for the solution of $\dcov \phi=0$.

We have the following:
\begin{itemize}
 \item {\textbf{$O(\lambda^{0})$}: The equation is $d\phi_{0}=0$, which by star-shapedness of $U$ is solved by
 \begin{equation}
    \phi_{0}=d\alpha_{0}
 \end{equation}
 for $\alpha_{0}\in \Lambda^{k-1}(U,V)$.}
 \item{\textbf{$O(\lambda^{1})$}: The equation is $d\phi_{1}+A\wedge \phi_{0}=0$. Taking $d$ on both sides of it we get $d(A\wedge \phi_{0})=0$, i.e., $A\wedge \phi_{0}=dH(A\wedge\phi_{0})$. We assume that this constraint fulfilled for now. That means, $d(\phi_{1}+H(A\wedge \phi_{0}))=0$, and the solution is
 \begin{equation}
  \phi_{1} = d\alpha_{1}-H(A\wedge\phi_{0}),
 \end{equation}
 for $\alpha_{1} \in \Lambda^{k-1}(U,V)$.
 }
 \item{\textbf{$O(\lambda^{2})$}: The equation is $d\phi_{2}+A\wedge\phi_{1}=0$, and the same procedure, with similar assumption, gives that the solution is
 \begin{equation}
  \phi_{2}=d\alpha_{2}-H(A\wedge\phi_{1}),
 \end{equation}
  for $\alpha_{2}\in \Lambda^{k-1}(U,V)$.}
  \item{\textbf{$O(\lambda^{l})$}: In general case the equation is $d\phi_{l}+A\wedge\phi_{l-1}=0$ which gives
  \begin{equation}
   \phi_{l}=d\alpha_{l}-H(A\wedge\phi_{l-1}),
  \end{equation}
  for $\alpha_{l}\in \Lambda^{k-1}(U,V)$.}
\end{itemize}
  Collecting all the terms we get for $\lambda=1$ the formal solution in terms of the series
  \begin{equation}
   \phi = (1-H(A\wedge\_)+H(A\wedge H(A\wedge\_))-\ldots) \sum_{l=0}^{\infty} d\alpha_{l}.
  \end{equation}
    Selecting $\{\alpha_{i}\}_{i=0}^{\infty}$ in such a way that their sum forms uniformly convergent series to a continuous function, denoting its sum by $\alpha:=\sum_{l=0}^{\infty} \alpha_{l}$, and taking into account that $dH\phi=d\alpha=c$ (using $H^{2}=0$), we get (\ref{Eq.Solution_homogenous_k_gt_0}). We can restrict $c$ to being smooth. If, in addition, $c \in ker(A\wedge\_)$ then we would get a contradiction with the assumption $\phi \not\in ker(A\wedge\_)$.

    %%%
    \textbf{Formal solution - integral equations approach:}
    %%%

    We present the other derivation of the solution (\ref{Eq.Solution_homogenous_k_gt_0}) using integral equation. If $\phi$ is a solution of (\ref{Eq_homogenous_equation}) then $A\wedge \phi = dH(A\wedge \phi)$, and we get that $d(\phi+H(A\wedge\phi))=0$, which results in the integral equation
    \begin{equation}
     \phi = - H(A\wedge\phi) + c,
    \end{equation}
    where $c$ is an arbitrary exact form. We can replace it with iteration scheme $\phi_{0}=c$ and $\phi_{i} = - H(A\wedge\phi_{i-1}) + c$, for $i > 0$, which gives the series that the formal limit $\lim_{i\rightarrow\infty}\phi_{i}$ is (\ref{Eq.Solution_homogenous_k_gt_0}).

    Taking the exterior derivative of each term of (\ref{Eq.Solution_homogenous_k_gt_0}) and summing formally, we get that
    \begin{equation}
     d\phi + A\wedge\phi = H(F\wedge \phi)=0,
    \end{equation}
    as required.

    %%%%
    \textbf{Convergence:}
    %%%%

    In the second part we prove that the formal series (\ref{Eq.Solution_homogenous_k_gt_0}) is convergent in some restricted neighbourhood of the homotopy center $x_{0}$. We estimate
    \begin{equation}
    \begin{array}{c}
     ||H(A\wedge \omega ) ||_{\infty} = || \int_{0}^{1}i_{\mathcal{K}}  (A\wedge\omega)(x_{0}+t(x-x_{0})) t^{k-1} dt ||_{\infty} \leq \int_{0}^{1}||x-x_{0}|| ||A||_{\infty} ||\omega||_{\infty}t^{k-1}dt \\
     = ||x-x_{0}|| ||A||_{\infty}||\omega||_{\infty}\frac{1}{k}.
    \end{array}
    \end{equation}
     We therefore have
    \begin{equation}
    \begin{array}{c}
     ||\phi||_{\infty} = ||(1- H(A\wedge\_)+H(A\wedge(H(A\wedge\_)))- \ldots)c||_{\infty} \leq  \\
     \left(1+||x-x_{0}|| \frac{||A||_{\infty}}{k} + \left(||x-x_{0}|| \frac{||A||_{\infty}}{k}\right)^{2} +\ldots \right) ||c||_{\infty},
    \end{array}
    \end{equation}
    and the series (\ref{Eq.Solution_homogenous_k_gt_0}) is absolutely and uniformly convergent to a continuous form when $||x-x_{0}||\frac{||A||_{\infty}}{k} <1$, which gives the condition (\ref{Eq_convergence_homogenous_solution}) for the radius of ball. Since $d\phi = -A\wedge \phi$ by the above construction, we have that it is a solution. Moreover, $d\phi$ is continuous.

    If $\phi\in ker(A\wedge\_)$ then from the equation $d\phi=0$ and so $\phi$ is a gauge mode. From (\ref{Eq.Solution_homogenous_k_gt_0}) we get that $\phi=c+\mathcal{A}(U,V)$, and therefor, $c \in \mathcal{E}(U,V)\cap ker(A\wedge\_)$ is also an exact gauge mode. Note that excluding exact gauge modes is essential since we would not obtain uniqueness.

    Moreover, from the preceding Proposition \ref{Proposition_noantiexactGaugeModes} we get that the only antiexact (and also exact) gauge mode is $\phi=0$.

    Uniqueness is proved in a standard way by reductio ad absurdum. Assume that there are two distinct solutions $\phi_{1} \neq \phi_{2}$ with the same initial conditions $dH\phi_{1}=dH\phi_{2}$. Then the form $\psi:=\phi_{1}-\phi_{2}$ is also the solution with $dH\psi=0$, i.e., the solution with no exact part. However, from the form of the solution obtained above, we see that if $dH\psi=0$, then $\psi=0$ (since $c=0$), so $\phi_{1}=\phi_{2}$, a contradiction.
\end{Proof}
%%%%%%%%%%%%%%%%%%%%%%%%%%%%%%%%%%
Note also the power series solution (\ref{Eq.Solution_homogenous_k_gt_0}) also works for $k=0$, when we set $c \in \mathbb{R}$ as an arbitrary constant.

\begin{Remark}
The solution is only continuous with exterior derivative continuous. This is what is expected when uniformly convergent series of smooth forms is used in constructing solutions. The similar situations is present in Picard-Lindel\"{o}f theorem for ODE, when the solution is $C^{1}$ class only\footnote{In our case the solution is continuous with exterior derivative continuous, which corresponds to continuous derivative for ODEs.}. Increasing smoothness imposes additional conditions, and is presented in Appendix \ref{Section_RegularityConsiderations}.

We can also assume that $A$ and $c$ is only of class $C^{1}$ without modifying the proof.
\end{Remark}

Moreover, one can note that the radius of convergence of (\ref{Eq.Solution_homogenous_k_gt_0}) is larger when the supremum of coefficients of $A$ is small.

From the above statements we have the following:

%%%%%%%%%%%%%%%%%%%%%%%%%%%%%%%%%%
\begin{Corollary}
The space of solution of $\dcov \phi =0$ is an affine space $\Lambda(U,V)\setminus ker(A\wedge\_)$ over the vector space of exact gauge modes $\mathcal{E}(U,V)\cap ker(A\wedge\_)$.
\end{Corollary}

%%%%%%%%%%%%%%%%%%%%%%%%%%%%%%%%%%
\begin{Remark}
 For each term of the series (\ref{Eq.Solution_homogenous_k_gt_0}) we have
 \begin{equation}
  A\wedge (H(A\wedge\_)^{l} c \in \mathcal{E}(U,V).
 \end{equation}
\end{Remark}
%%%%%%%%%%%%%%%%%%%%%%%%%%%%%%%%%%
Moreover,
%%%%%%%%%%%%%%%%%%%%%%%%%%%%%%%%%%
\begin{Remark}
We can formally write (\ref{Eq.Solution_homogenous_k_gt_0}) as
\begin{equation}
 \phi = \frac{1}{1+HA\wedge\_}c.
 \label{Eq_inverison_of_1HA}
\end{equation}
This notation will be firmly stated within the framework of Operational Calculus later.
\end{Remark}
%%%%%%%%%%%%%%%%%%%%%%%%%%%%%%%%%%
%%%%%%%%%%%%%%%%%%%%%%%%%%%%%%%%%%
\begin{Remark}
 In solving (\ref{Eq_homogenous_equation}) the 'initial condition' for the iterative procedure described by the series (\ref{Eq.Solution_homogenous_k_gt_0}) is a form $c \in \mathcal{E}(U,V)) \setminus ker(A\wedge\_)$. In this sense the exact form $c$ parametrizes the solution.
\end{Remark}
%%%%%%%%%%%%%%%%%%%%%%%%%%%%%%%%%%

We can now formulate the algorithm for solving (\ref{Eq_homogenous_equation}).
%%%%%%%%%%%%%%%%%%%%%%%%%%%%%%%%%%
\begin{Algorithm}
\label{Algorithm_1}
 In order to solve
 \begin{equation}
  \dcov \phi =0,
 \end{equation}
 for $\phi\in \Lambda^{k}(U,V)\setminus ker(A\wedge\_)$, $k>0$, pick an initial condition $\gamma_{0}\in \mathcal{E}(U, V)\setminus \ker(A\wedge\_)$, and compute iteratively
 \begin{equation}
  \gamma_{l}=H(A\wedge \gamma_{l-1}).
 \end{equation}
 Then the solution is
 \begin{equation}
  \phi = \sum_{l=0}^{\infty} (-1)^{l} \gamma_{l}.
 \end{equation}
The series is convergent for such $x \in B(x_{0},r)$, where  $r\frac{||A||_{\infty}}{k}=1$.

To this solution you can add an arbitrary exact gauge mode from $\mathcal{E}(U,V)\cap ker(A\wedge\_)$.
\end{Algorithm}
%%%%%%%%%%%%%%%%%%%%%%%%%%%%%%%%%%

We now provide a simple example that explains the Algorithm \ref{Algorithm_1}.
%%%%%%%%%%%%%%%%%%%%%%%%%%%%%%%%%%
\begin{Example}
\label{Ex1}
 Let us solve the equation
 \begin{equation}
  d\phi+A\wedge\psi=0,
 \end{equation}
 on $\mathbb{R}^{2}$ with coordinates $(x,y)$, where $A =dy$, and with initial condition $\gamma_{0}=dx\in \mathcal{E}(\mathbb{R}^2)$ and the center $x_{0}=0$, i.e., $\mathcal{K}=x\partial_{x}+y\partial_{y}$. Note that $\gamma_{0}\not\in ker(A\wedge\_)$.

 We have
 \begin{itemize}
  \item {$\gamma_{1} = \int_{0}^{1}\mathcal{K} \lrcorner \left( dy\wedge dx \right) t dt = \frac{1}{2!}(ydx-xdy)$,}
  \item {$\gamma_{2} = \int_{0}^{1}\mathcal{K} \lrcorner \left(\frac{1}{2!}y dy\wedge dx \right) t^{2} dt  = \frac{1}{3!}(y^{2}dx-yxdy)$,}
  \item {$\gamma_{3} = \int_{0}^{1} \mathcal{K}\lrcorner \left( \frac{1}{3!} y^{2} dy\wedge dx \right) t^{3} dt = \frac{1}{4!}(y^{3}dx-y^{2}xdy)$,}
  \item {$\gamma_{k}=\frac{1}{(k+1)!}\left(y^{k}dx-y^{k-1}xdy \right)$.}
 \end{itemize}
Then, by summing the terms, we get
\begin{equation}
 \phi = \sum_{l=0}^{\infty} \gamma_{l} = (1- e^{-y})\frac{dx}{y} + (e^{-y}-1+y)\frac{xdy}{y^{2}}.
 \label{Ex1.solution}
\end{equation}
The solution has a removable singularity at $y=0$.

The projection to the initial condition is given by $dH$, and we have $dH \phi = dx$ as required.

By straightforward computations we have
\begin{equation}
 \begin{array}{c}
  d\phi = \frac{1}{y}(1-e^{-y})dx\wedge dy, \\
  A\wedge \phi = \frac{1}{y}(1-e^{-y})dy\wedge dx,
 \end{array}
\end{equation}
so $\dcov \phi=0$ as required.

One can note that the solution $\phi$ is well-defined for the whole $\mathbb{R}^{2}$, so its radius of convergence is significantly larger than the Theorem \ref{Th_homogenous_solution} suggests.

Moreover, if we treat $\mathbb{R}^{2}$ as a fibered bundle with horizontal direction $dx$ and vertical $dy$,  then the form (\ref{Ex1.solution}) is neither horizontal nor vertical. So projecting $\phi$ to $dx$ component and then lifting back along $A=dy$ we do not obtain the original form and even their covariant constancy.

One can also note that if we would choose $\gamma_{0}=f(x)dy \in ker(A\wedge\_)$ such that $\gamma_{0}\not\in \mathcal{E}(\mathbb{R}^{2})$, i.e., $\partial_{x}f(x) \neq 0$, then we would get a contradiction.

We can also note that the gauge mode has the form $f(y)dy$, for arbitrary $f \in \Lambda^{0}(\mathbb{R}^2)$, and we can add it to the solution.

\end{Example}
%%%%%%%%%%%%%%%%%%%%%%%%%%%%%%%%%%
\begin{Example}
\label{Ex2}
 Continuing Example \ref{Ex1}, we can also check easily (by assuming the solution in the form $\phi_{2}=f(y)dx$) that the solution is
 \begin{equation}
  \phi_{2}=e^{-y}dx.
 \end{equation}
 It is interesting to note that the initial condition for this solution is the solution from the previous example since
 \begin{equation}
  dH\phi_{2}=\phi = (1- e^{-y})\frac{dx}{y} + (e^{-y}-1+y)\frac{xdy}{y^{2}}.
 \end{equation}
 This is the condition for starting the algorithm to obtain $\phi_{2}$.

The solution $\phi_{2}$ is horizontal when treating $dx$ as a horizontal direction.
\end{Example}

%%%%%%%%%%%%%%%%%%%%%%%%%%%%
Finally, we provide an example with a non-scalar connection.
\begin{Example}
Consider the covariant derivative of the form
\begin{equation}
 \dcov =d+A\wedge\_ = d+\left[ \begin{array}{cc}
           \alpha & 0 \\
           0 & 0 \\
          \end{array}
\right]\wedge\_,
\end{equation}
for $\alpha \in \Lambda^{1}(U)$.

Consider the case $\Lambda^{k}(U,V)$ with $k>0$. The kernel of $A\wedge\_$ consists of a set of vectors
\begin{equation}
 \phi=\left[\begin{array}{c}
        \beta \\
        \gamma
       \end{array}
       \right],
\end{equation}
where $\beta\in \{ \psi | \alpha\wedge \psi=0\}$, and $\gamma$ is an arbitrary differential form.

Excluding these forms we can uniquely solve the equation $\dcov \phi=0$ using (\ref{Eq.Solution_homogenous_k_gt_0}) for initial data $c \not\in \mathcal{E}(U,V)\cap ker(A\wedge\_)$, $c\neq0$. To such a solution we can add an arbitrary exact gauge mode.
\end{Example}

%%%%%%%%%%%%%%%%%%%%%%%%%%%%%%%%%%

%%%%%%%%%%%%%%%%
\subsection{Inhomogenous equation}
\label{Subsection_Inhomogenous_equation}
%%%%%%%%%%%%%%%%

The next step is to provide a solution for the inhomogeneous covariant constancy equation. An intuitive discussion of the influence of inhomogeneous terms is given in Appendix \ref{Appendix_Antiexact_inhomogenity}.

%%%%%%%%%%%%%%%%%%%%%%%%%%%%
\subsection{Exact inhomogenity}
%%%%%%%%%%%%%%%%%%%%%%%%%%%%

We begin with the particular case when the inhomogeneity is an exact form, i.e., for equation
\begin{equation}
 \dcov \phi = J_{e},
 \label{Eq_inhomogenousExactequation}
\end{equation}
where $J_{e}\in\mathcal{E}(U,V)$.

We note that:
\begin{Remark}
When $\phi \in \Lambda(U,V)$ is a solution of (\ref{Eq_inhomogenousExactequation}) then we can decompose $\phi=\phi_{1}+\phi_{2}$, where $\phi_{2}\in ker(A\wedge\_)$, i.e. solves
\begin{equation}
 A\wedge\phi_{2} =0,
\end{equation}
and $\phi_{1}\in\Lambda(U,V)\setminus ker(A\wedge\_)$ solves
\begin{equation}
 d\phi_{1} +A\wedge \phi_{1} = J_{e} -d\phi_{2}.
\end{equation}
Therefore if $\phi_{2}\in\mathcal{E}(U,V)\cap ker(A\wedge\_)$ then it is a gauge mode that does not modify inhomogeneity $J_{e}$ since $d\phi_{2}=0$. If however $\phi_{2}\in\mathcal{A}(U,V)\cap ker(A\wedge\_)$, then $d\phi_{2}$ modifies inhomogenity since then $d\phi_{2}\neq 0$.
\end{Remark}

We first consider the scalar case:
%%%%%%
\begin{Proposition}
 The unique smooth solution $\phi\in \Lambda^{0,\infty}(U,\mathbb{R})$ of
 \begin{equation}
  \dcov \phi  = J_{e},
  \label{Eq_nonhomogenous_covariant_equation_exact}
 \end{equation}
for $A\in \Lambda^{1, \infty}(U, End(\mathbb{R}))$, $J_{e}\in \mathcal{E}^{1, \infty}(U)$, with $c\in \mathbb{R}$ is
\begin{equation}
 \phi = \exp(-HA)(c+H(J_{e}\exp(HA)).
 \label{Eq_solution_nonhomogenous_k_eq_0}
\end{equation}
\end{Proposition}
\begin{Proof}
 We have the solution of the homogenous equation $\phi=C\exp(-HA)$. By the variation of the constant, i.e., taking $C\in \Lambda^{0, \infty}(U)$, and substituting back to the equation, we obtain
 $dC = J_{e}\exp(HA)\in \mathcal{E}(U)$, and as a result, $J_{e}\exp(HA)=dH(J_{e}\exp(HA))$. This gives $d(C-H(J_{e}\exp(HA)))=0$, i.e, $C=D+H(J_{e}\exp(HA))$, for real number $D$. This, replacing $D$ with $c$, gives (\ref{Eq_solution_nonhomogenous_k_eq_0}).
\end{Proof}

We split the investigation into two cases. First one is when $\phi \in ker(A\wedge\_)$.
%%%%%%%%%%%%%%%%%%%%%%%
\begin{Proposition}
 The solution of $\dcov \phi = J_{e}$ for $J_{e}\in \mathcal{E}^{k+1}(U,V)$, $k>0$ is an element of $ker(A\wedge\_)$ when
 \begin{equation}
 \phi = c + HJ_{e},
 \end{equation}
with $c \in \mathcal{E}(U,V)\cap ker(A\wedge\_)$ is an exact gauge mode, and with $HJ_{e}\in  \mathcal{A}(U,V)\cap ker(A\wedge\_)$ is an antiexact element. The second condition is the constraint on $J_{e}$.
\end{Proposition}
\begin{Proof}
If $\phi\in ker(A\wedge\_)$, then it must fulfill $d\phi = J_{e}$. Since $J_{e}$ is exact ($J_{e}=dHJ_{e})$, so $\phi = c + HJ_{e}$ for some exact $c$. Since the condition $A\wedge \phi=0$ is linear, so both $c$ and $HJ_{e}$ must belong also to $ker(A\wedge\_)$.
\end{Proof}

Now we provide the solution for $\phi \not\in ker(A\wedge\_)$.
%%%%%%%%%%%%%%%%%%%%%%
\begin{Theorem}
\label{Th_nonhomogenous_solution_exactRHS}
 The unique solution $\phi \in \Lambda^{k}(U, V)\setminus ker(A\wedge\_)$ of
 \begin{equation}
  \dcov \phi  = J_{e},
  \label{Eq_nonhomogenous_covariant_equation}
 \end{equation}
for $A\in \Lambda^{1}(U, End(V))$, $J_{e}\in \mathcal{E}^{k+1}(U,V)$ for $k>0$, with $dH\phi=c\in \mathcal{E}(U,V)\setminus ker(A\wedge\_)$ is
\begin{equation}
 \phi=\phi_{H} + \phi_{I}, \quad \phi_{I}=\sum_{l=0}^{\infty}(-1)^{l} (H(A\wedge\_))^{l} HJ_{e},
 \label{Eq_solution_nonhomogenous_k_gt_0}
\end{equation}
where $\phi_{H}$ is a solution of homogenous equation ($J_{e}=0$) given in Theorem \ref{Th_homogenous_solution}.

The series in (\ref{Eq_solution_nonhomogenous_k_gt_0}) is convergent to the continuous solution in the ball $B(x_{0},r)$, where the radius is given by $r \frac{||A||_{\infty}}{k}=1,$, where the supremum norm is taken over the closed ball.

When $HJ_{e}\in \mathcal{E}^{k+1}\cap ker(A\wedge\_)$, then the solution reduces to
\begin{equation}
 \phi = \phi_{H}+HJ_{e}.
\end{equation}
\end{Theorem}
%%%%%%%%%%%%%%%%%%%%%%%%%%%%%%%%%%

%%%%%%%%%%%%%%%%%%%%%%%%%%%%%%%%%%
\begin{Proof}
We proceed as in the proof of the previous theorem.

In the first part we derive a formal solution using formal perturbation series for differential equations. We replace the equation (\ref{Eq_nonhomogenous_covariant_equation}) by
 \begin{equation}
  d\phi+\lambda A\wedge\phi = J_{e},
 \end{equation}
 for a nonzero real number $\lambda$. Introducing the formal ansatz $\phi=\sum_{l=0}^{\infty}\lambda^{l}\phi_{l}$ we get:
 \begin{itemize}
  \item {\textbf{$O(\lambda^{0})$}: The equation is $d\phi_{0}=J\in \mathcal{E}$, and therefore, $J_{e}=dHJ_{e}$, so
  \begin{equation}
   \phi_{0}=d\alpha_{0}+HJ_{e},
  \end{equation}
    for $\alpha_{0}\in \Lambda^{k-1}(U,V)$. }
 \item { \textbf{$O(\lambda^{l})$}: We get the recurrence for the solution
 \begin{equation}
  \phi_{l}=d\alpha_{l} - H(A\wedge \phi_{l-1}),
 \end{equation}
    for $\alpha_{l}\in \Lambda^{k-1}(U,V)$.}
 \end{itemize}
  Summing up the terms we have
  \begin{equation}
   \phi = \underbrace{ \sum_{l=0}^{\infty}(-1)^{l} (H(A\wedge\_))^{l} \sum_{p=0}^{\infty}d\alpha_{p}  }_{\phi_{H}} + \underbrace{\sum_{l=0}^{\infty}(-1)^{l}(H(A\wedge\_))^{l} HJ_{e}}_{\phi_{I}}.
  \end{equation}
    As before, we can select $\{\alpha_{p}\}_{p=0}^{\infty}$ to form uniformly convergent series, so setting $\alpha=\sum_{p=0}^{\infty}\alpha_{p}$ and using the condition $dH\phi=d\alpha=c$ we get (\ref{Eq_solution_nonhomogenous_k_gt_0}).

  The second part is the proof of convergence and uniqueness of the formal series. It is the same as in the proof for the homogenous case.

  One can also note that $\phi \not\in ker(A\wedge\_)$ since otherwise it would be that $c \in ker(A\wedge\_)$ (as well as $HJ_{e}$). Then we would get nonuniqueness since $\phi_{H}$ would be nonunique.
\end{Proof}
%%%%%%%%%%%%%%%%%%%%%%%%%%%%%%%%%%
From the proof we have
%%%%%%%%%%%%%%%%%%%%%%%%%%%%%%%%%%
\begin{Remark}
 The solution (\ref{Eq_solution_nonhomogenous_k_gt_0}) can be written as
 \begin{equation}
  \phi=\phi_{H}+G(J_{e}),
 \end{equation}
 where $G$ resembles a Green's function used in the theory of the second order Laplace-Beltrami operator $\triangle = d\delta+\delta d$, see, e.g., \cite{Thirring}. However, we do not assume a metric structure here, so the approach is general. Moreover, no boundary conditions were imposed on $G$.
\end{Remark}
%%%%%%%%%%%%%%%%%%%%%%%%%%%%%%%%%%
Finally, we can consider the most general case.
%%%%%%%%%%%%%%%%%%%%%%%%%%%%%%%%%%
\begin{Theorem}
The most general solution of
\begin{equation}
 \dcov \phi = J_{e},
\end{equation}
for $\phi \in \Lambda^{k}(U,V)$, $J_{e}\in \mathcal{E}^{k+1}(U,V)$ is of the form
\begin{equation}
\phi = \phi_{1}+\phi_{2},
\end{equation}
where $\phi_{1} \not \in ker(A\wedge\_)$ is the solution of the inhomogeneous equation with exact RHS (see Theorem \ref{Th_nonhomogenous_solution_exactRHS})
\begin{equation}
 d\phi_{1}+A\wedge\phi_{1} =J_{e}-d\phi_{2},
\end{equation}
and $\phi_{2}\in ker(A\wedge\_)$ is the solution of
\begin{equation}
 A\wedge \phi_{2}=0.
\end{equation}
\end{Theorem}
%%%%%%%%%%%%%%%%%
\begin{Proof}
 Since the operator $A\wedge\_$ is linear, we can decompose $\phi=\phi_{1}+\phi_{2}+\phi_{3}$, where $\phi_{1}, \phi_{3} \not\in ker(A\wedge\_)$, $\phi_{2}\in ker(A\wedge\_)$ in such a way that
 \begin{equation}
  A\wedge\phi_{1}\in \mathcal{E}^{k+1}(U,V),~ A\wedge\phi_{3}\in \mathcal{A}^{k+1}(U,V).
 \end{equation}
 Substituting back to the equation, we have that
 \begin{equation}
 \underbrace{d\phi_{1}+A\wedge\phi_{1}+d\phi_{2}-J_{e} +d\phi_{3}}_{\mathcal{E}} + \underbrace{A\wedge\phi_{3}}_{\mathcal{A}}=0,
 \end{equation}
where exact ($\mathcal{E}$) and antiexact ($\mathcal{A}$) parts were marked. Both direct summands must vanish, so we get that $\phi_{3}=0$ and finish the proof.
\end{Proof}

Note that since we do not require that $\phi \not\in ker(A\wedge\_)$ so, we do not expect the uniqueness, since the exact gauge modes as well as elements of $ker(A\wedge\_)$ can be present.

We will provide examples.
%%%%%%%%%%%%%%%%%%%%%%%%%%%%%%%%%%
\begin{Example}
 Continuing Example \ref{Ex1}, we solve the equation
 \begin{equation}
  \dcov \phi = J, \quad J=xdx, \quad A=dy,
 \end{equation}
where $J$ is an exact form. We have
\begin{equation}
 HJ = \frac{1}{2}x^{2}.
\end{equation}

First, we use the series solution (\ref{Eq_solution_nonhomogenous_k_gt_0}) and then we compare it with (\ref{Eq_solution_nonhomogenous_k_eq_0}). We have
\begin{itemize}
 \item {$HA\wedge HJ=H(\frac{1}{2}x^{2}dy)=\frac{1}{2}x^{2}y\int_{0}^{1}t^{2}dt = \frac{1}{3!}x^{2}y$,}
 \item {$(HA\wedge)^{2}HJ=\frac{1}{4!}x^{2}y^{2}$,}
 \item {$\ldots$.}
 \item {$(HA\wedge)^{k}HJ=\frac{1}{(k+2)!}x^{2}y^{k}$.}
\end{itemize}
The inhomogenous part of the solution (\ref{Eq_solution_nonhomogenous_k_gt_0}) is given by the series
\begin{equation}
 \phi_{I}=\frac{1}{2!}x^{2}-\frac{1}{3!}x^{2}y+\ldots = \left(\frac{x}{y}\right)^{2}(e^{-y}-1+y).
\end{equation}

Likewise, applying (\ref{Eq_solution_nonhomogenous_k_eq_0}), we have
\begin{equation}
\begin{array}{c}
 \phi_{I}=\exp(-HA) H(J\exp(HA)) = e^{-y}H(e^{y}xdx)=e^{-y}x^{2}\int_{0}^{1}te^{ty}dt = \\
 x^{2}\frac{d}{dy}\frac{1}{y}\int_{0}^{y}e^{z}dz = \left(\frac{x}{y}\right)^{2}(e^{-y}-1+y),
\end{array}
\end{equation}
as previously. Therefore the inhomogeneous contributions calculated either by (\ref{Eq_solution_nonhomogenous_k_gt_0}) or (\ref{Eq_solution_nonhomogenous_k_eq_0}) agree, as required.
\end{Example}

%%%%%%%%%%%%%%%%%%%%%%%%%%%%%%%%%%%%%%%%%%%%
\subsection{General inhomogenity}
%%%%%%%%%%%%%%%%%%%%%%%%%%%%%%%%%%%%%%%%%%%%

Finally, we will consider the most general equation where $J$ is an arbitrary, not necessarily exact, form. We have
%%%%%%%%%%%%%%%%
\begin{Theorem}
\label{Th_FullInhomogenous_parallelTransportEquation}
 The solution of the inhomogeneous covariant constancy equation
 \begin{equation}
  \dcov \phi = J,\quad \dcov = d + A\wedge \_,
  \label{Eq.FullInhomogenous_CovariancyConstantEquation}
 \end{equation}
where $\phi\in\Lambda^{k}(U, V)$, $A\in \Lambda^{1}(U,End(V))$, $J\in \Lambda^{k+1}(U,V)$, $k>0$ is given by
\begin{equation}
 \phi = \phi_{1}+\phi_{2}+\phi_{3},
\end{equation}
where $\phi_{1}$ fulfils
\begin{equation}
 \dcov \phi_{1}=J_{e} - d(\phi_{2}+\phi_{3}),
\end{equation}
$\phi_{2}$ fulfils
\begin{equation}
 A\wedge \phi_{2} = J_{a},
 \label{Eq_Ja_condition}
\end{equation}
 where $J_{e}:=dHJ$ is the exact part of $J$, and $J_{a}:=HdJ$ is the antiexact part of $J$.
 Finally, we have an arbitrary choice for $\phi_{3}\in \ker(A\wedge\_)$. Moreover, $A\wedge \phi_{1}\in \mathcal{E}^{k+1}(U, V)$ and $A\wedge\phi_{2}\in \mathcal{A}^{k+1}(U,V)$.

 If the equations cannot be solved, then there is no solution of (\ref{Eq.FullInhomogenous_CovariancyConstantEquation}). The only equation that imposes constraint is (\ref{Eq_Ja_condition}) and can be fulfilled only when $J_{a}\in Im(A\wedge\_)$.
\end{Theorem}
%%%
\begin{Proof}
 Since $A\wedge \phi \in \Lambda^{l}(U,V)$ and $l>0$, so it can be decomposed into exact and antiexact parts. Therefore, we can find three forms $\phi=\phi_{1}+\phi_{2}+\phi_{3}$ such that
 \begin{equation}
  A\wedge \phi_{1} \in \mathcal{E}(U,V), \quad A\wedge \phi_{2}\in \mathcal{A}(U,V), \quad A\wedge \phi_{3}=0.
 \end{equation}
Decomposing $J=J_{e}+J_{a}$ and substituting into the equation (\ref{Eq.FullInhomogenous_CovariancyConstantEquation}) and splitting into exact ($\mathcal{E}$) and antiexact ($\mathcal{A}$) parts gives
\begin{equation}
 \underbrace{d(\phi_{1}+\phi_{2}+\phi_{3})+A\wedge\phi_{1} -J_{e}}_{\mathcal{E}} +\underbrace{A\wedge \phi_{2} - J_{a}}_{\mathcal{A}}=0.
\end{equation}
Using direct sum decomposition into exact and antiexact terms ends the proof.
\end{Proof}
%%%%%%%%%%%%
\begin{Remark}
We can note that the existence of $J_{a}$ can be a severe obstruction to the existence of the solution to the general equation. If $J_{a}=0$ then the solution (not necessarily unique) exists. The extraction of $J_{a}$ and examination of $Im(A\wedge\_)$, therefore, is a test for the existence of a solution for a general case.
\end{Remark}

\begin{Remark}
 If exists, $\phi_{2}=\phi_{2}(J_{a})$, i.e., $\phi_{2}$ depends on the antiexact part of $J_{a}$.
\end{Remark}

We can now construct a practical way of solving the equation (\ref{Eq.FullInhomogenous_CovariancyConstantEquation}).
%%%%%%%%%%%%%%%%%
\begin{Algorithm}
 For the equation (\ref{Eq.FullInhomogenous_CovariancyConstantEquation}):
 \begin{enumerate}
     \item {determine kernel of $A\wedge\_$:
    \begin{equation}
     A\wedge \phi_{3}=0,
    \end{equation}
    for $\phi_{3}$ and pick one element,
    }
  \item {solve the algebraic constraint:
  \begin{equation}
   A\wedge \phi_{2} = J_{a},
  \end{equation}
    for $\phi_{2}\not \in ker(A\wedge\_)$; If the solution does not exist then there is no solution of (\ref{Eq.FullInhomogenous_CovariancyConstantEquation});}
    \item {solve the differential equation:
    \begin{equation}
     d\phi_{1}+A\wedge\phi_{1}=J_{e}-d(\phi_{2}+\phi_{3}),
    \end{equation}
    for $\phi_{1}$ which is an inhomogenous covariant constant equation with exact RHS,
    }
    \item {
    compose the full solution
    \begin{equation}
     \phi=\phi_{1}+\phi_{2}+\phi_{3}.
    \end{equation}
    }
 \end{enumerate}
\end{Algorithm}
%%%%%%%%%%%%%%%%%%%%%%

We provide a simple example of solving algebraic constraints. First, the negative example will be provided
\begin{Example}
 We continue the Example \ref{Ex1}. We consider
 \begin{equation}
  \dcov \phi=J_{a}, \quad J_{a}=\frac{1}{2}(xdy-ydx), \quad A=dy.
 \end{equation}
 Since the solution of exact part, from Theorem \ref{Th_homogenous_solution}, is $\phi_{1}=ce^{-y}$ and $\phi_{3}=0$, we will focus only on algebraic constraint
 \begin{equation}
  A\wedge \phi_{2} = J_{a}.
  \label{Ex_Jaconstraint_0}
 \end{equation}
 Assuming that $\phi_{2}=f(x,y) \in \Lambda^{0}(U)$, substituting into (\ref{Ex_Jaconstraint_0}) we get a contradiction. Therefore, there are no solutions to this problem.
\end{Example}
%%%%%%%%%%%%%%%%%%%%
As a positive example, we propose the following
\begin{Example}
 Consider a star-shaped $U \subset \mathbb{R}^{3}$ with coordinates $x,y,z$. We will try to solve the constraint
 \begin{equation}
  A\wedge \phi = J_{a}, \quad J_{a}=xdy\wedge dz -y dx\wedge dz+z dx\wedge dy, \quad A=dy,
 \end{equation}
 which is part of solving the full equation $\dcov \phi=J_{a}$. To solve this constraint, assume that
\begin{equation}
 \phi =f(x,y,z)dx + g(x,y,z)dy + h(x,y,z) dz,
\end{equation}
for $f,g,h \in \Lambda^{0}(U)$. By substituting into constraints, one gets
\begin{equation}
 f = -z, \quad h=x,
\end{equation}
and $g$ is arbitrary. We set $g=0$ to remove $\phi$ that are elements of $ker(A\wedge\_)$.
\end{Example}
%%%%%%%%%%%%%%

%%%%%%%%%%%%%%%%%%%%%%%%%%%%%%%%%%%%%%%%%%%%%%%%%%%%

As a final remark of this section, note that the operators for the inverse of $\dcov$ are nonlocal. They can be expressed by curvature using
%%%%%%%%%%%%%%%%%%%%%%%%%%%%%%%%%%%%%%%%%%%%%%%%%%%%
\begin{Proposition}
\begin{equation}
 H(A\wedge d\alpha) = dH(A\wedge \alpha) + H(F\wedge\alpha)-H(A\wedge A\wedge\alpha)  -A\wedge \alpha
\end{equation}
where $F=\dcov A$ is the curvature. For the connection valued in abelian groups $A\wedge A=0$.
\end{Proposition}
%%%%%%%%%%%%%%%%%%%%%%%%%%%%%%%%%%
Using the above proposition and making a recursive substitution 'inside-out' in (\ref{Eq.Solution_homogenous_k_gt_0}), we are led to a complicated expression, and therefore, we do not follow this path. One can notice that the curvature (and connection form) enter the solution in a highly nonlocal and nonlinear way.

In the next section we discuss the issue of basic forms on associated vector bundles.

%%%%%%%%%%%%%%%%%%%%%%%%%%%%%%%%%%%%%%%%%%%%%%%%%%%%%%%%%%%%%%%%%
\section{Associated vector bundles}
%%%%%%%%%%%%%%%%%%%%%%%%%%%%%%%%%%%%%%%%%%%%%%%%%%%%%%%%%%%%%%%%%
On the associated vector bundle, the base (horizontal and equivariant) vector-valued forms are in 1:1 correspondence with sections of the bundle. Therefore, we must establish the correspondence between solutions from our method with horizontal projection and equivariance. We present it in the following two subsections.

%%%%%%%%%%%%%%%%%%%%%%%%%%%%%%%%%%%%%%%%%%%%%%%%%%%%%%%%%%%%%%%%%
\subsection{Horizontal projection}
%%%%%%%%%%%%%%%%%%%%%%%%%%%%%%%%%%%%%%%%%%%%%%%%%%%%%%%%%%%%%%%%%
We analyze possible problems when considering the horizontality and covariance constancy of forms on fibered sets. This issue is essential if we want to relate the solutions on the fibered set/bundle to the forms on base space as in the case of the associated bundle. Since we only want to illustrate the issue, we consider only scalar-valued forms from $\Lambda(U)$ for simplicity. For vector-valued forms, additional constraints related to the matrix structure of the $A$ form must be considered. The idea of this section is based on an adaptation of the proof of the Retraction theorem (Theorem 1.3 of \cite{Bryant}).

Within the setup of fibered space $U$ and a one-form $A\in \Lambda^{1}(U)$ let us call $ TU\setminus ker(A)=VU$ the vertical tangent space with dimension $k=dim(VU)$.

\begin{Proposition}
 For a one-form $A$ we can select $k=dim(VU)$ linearly independent vectors $\{X_{i}\}_{i=1}^{k}$ such that
 \begin{equation}
  X_{i}\lrcorner A = 1.
  \label{Eq_Xi_definition}
 \end{equation}
 Moreover, the one-form $A$ can be decomposed into the sum of linearly independent one-forms
\begin{equation}
 A = \sum_{i=1}^{k}\omega_{i}, \quad \omega_{1}\wedge\ldots\wedge \omega_{k}\neq 0.
\end{equation}
For each such one-form $\omega_{i}$ we can select a vector $X_{i}$ such that
\begin{equation}
 X_{j}\lrcorner\omega_{i}=\delta_{ij}.
 \label{Eq_omega_orthogonality}
\end{equation}

\end{Proposition}
\begin{Proof}
For the proof assume that there is $X_{i}$ such that in addition we have $X_{i}\lrcorner \omega_{j}=a\neq 0$. Then since $\omega_{i}$ and $\omega_{j}$ are linearly independent, we get a contradiction. Linear independence of vectors can be proved similarly.
\end{Proof}

Now the vertical space $VU = span(\{X_{i}\}_{i=1}^{k})$. We can construct projectors:
\begin{equation}
 P_{i} = I-\omega_{i}\wedge (X_{i}\lrcorner \_),
\end{equation}
with the property
\begin{equation}
 X_{i}\lrcorner P_{i} =0.
\end{equation}

We can see that, since (\ref{Eq_omega_orthogonality}) is valid, we have that the operators $\{P_{i}\}$ commutes pairwise, i.e.,
\begin{equation}
 P_{i}\circ P_{j}=P_{j}\circ P_{i}.
\end{equation}

The projectors $P_{i}$ are homomorphism of exterior algebra since we have
\begin{equation}
 P_{i}(\alpha\wedge\beta)=P_{i}(\alpha)\wedge P_{i}(\beta),
\end{equation}
where $\omega_{i}\wedge \omega_{i}=0$ was used. It can also be noted that for $\omega_{i}$ (and therefore $A$) being non-scalar this property does not generally hold.

We can now project a differential form $\phi \in \Lambda^{*}(U)$ onto a horizontal form by means of
\begin{equation}
 \Delta:=P_{1}\circ \ldots \circ P_{k}.
\end{equation}
that is
\begin{equation}
 X_{i} \lrcorner \Delta \phi =0,
\end{equation}
for all $1<i<k$.

Now we can examine the relation between solutions of $\dcov \phi =0$ and its horizontal part $\Delta \phi$. We have an obvious statement that if $[\dcov, \Delta ] = 0$, then the horizontal part of $\phi$ is also covariantly constant. However, generally, this is not the case. To grasp more knowledge about the problem notice that we can write
\begin{equation}
 \Delta = I -\sum_{i} \omega_{i}\wedge X_{i}\lrcorner\_ - \sum_{i<j} \omega_{i}\wedge\omega_{j} \wedge (X_{i}\lrcorner X_{j}\lrcorner\_) + \sum_{i<j<l} \omega_{i}\wedge \omega_{j}\wedge\omega_{l}\wedge (X_{i}\lrcorner X_{j} \lrcorner X_{l}\lrcorner \_)+\ldots,
\end{equation}
where all summation indices run over the set $\{1,\ldots, k\}$. Therefore if $\dcov \phi=0$ then $\dcov \Delta \phi=0$ if
\begin{equation}
 \dcov \Delta\phi = \sum_{i} \dcov(\omega_{i}\wedge X_{i}\lrcorner\phi) + \sum_{i<j} \dcov(\omega_{i}\wedge\omega_{j}\wedge (X_{i}\lrcorner X_{j}\lrcorner \phi)) + \ldots = 0.
\end{equation}
By vanishing all the summands, we get the necessary conditions for the covariantly constant solutions to be horizontal.

We illustrate it using an example.
\begin{Example}
 Continuing Example \ref{Ex1}, we have $X_{1}=\partial_{y}$ and
 \begin{equation}
  \Delta \phi = (1-e^{-y})\frac{dx}{y}.
 \end{equation}
 However, $\dcov \Delta \phi \neq 0$.

In Example \ref{Ex2} we have that $\Delta \phi_{2}=\phi_{2}$, so $P_{1}|_{\phi_{2}} = I$ and therefore $\dcov \Delta \phi_{2}=\dcov \phi_{2}=0$. Therefore, in this case $\Delta$ operator commute with $\dcov$.
\end{Example}

%%%%%%%%%%%%%%%%%%%%%%%%%%%%%%%%%%%%%%%%%%%%%%%%%%%%%%%%%%%%%%%%%
\subsection{Equivariance}
%%%%%%%%%%%%%%%%%%%%%%%%%%%%%%%%%%%%%%%%%%%%%%%%%%%%%%%%%%%%%%%%%
Under the gauge transformation by $g \in G$ the connection one-from transforms as
\begin{equation}
 A' = Ad(g^{-1})A + g^{-1}dg.
\end{equation}
Denoting by $\dcov$ the exterior covariant derivative induced by $A$ and by $d^{\nabla '}$ this induced by $A'$, we have
\begin{Proposition}
If $\phi\in \Lambda(U,V)\setminus ker(A\wedge\_)$ is a solution of $\dcov$, then $\phi' = g \phi \Lambda(U,V)\setminus ker(A'\wedge\_)$ is a solution of $d^{\nabla '}$.
\end{Proposition}
\begin{Proof}
The proof idea is similar to proposition 3.3 of \cite{SuperPoincareLemma}: We show that RHS of $\phi' = g^{-1} \phi$ fulfills parallel transport equation, so by the uniqueness of the solution of the homogenous equation, cf. Theorem \ref{Th_homogenous_solution}, we have that both sides are equal. Note that the uniqueness can be used only when we exclude from solutions gauge modes: $ker(A\wedge\_)\cap \mathcal{E}(U,V)$, as discussed before.

By a standard argument, we have
\begin{equation}
 \begin{array}{c}
  d^{\nabla '}\phi' = d(g^{-1}\phi)+ \left( g^{-1} A g + g^{-1}dg\right)\wedge g^{-1} \phi = \\
  dg^{-1} \wedge \phi + g^{-1} d\phi + g^{-1}A\wedge \phi + g^{-1}dg g^{-1} \wedge \phi = \\
  g^{-1} \dcov \phi,
 \end{array}
\end{equation}
where $0 = d(g^{-1}g) = d(g^{-1})g + g^{-1}dg$ was used. Since $\dcov \phi =0$ so $d^{\nabla '}\phi' =0$.
\end{Proof}

Since gauge transformation is equivalent to the move between horizontal sections of a principal bundle of the associated bundle, we get the equivariance of $\phi$, i.e., $\phi(pg)=g^{-1}\phi(p)$.

The following section connects the results obtained so far with Operational Calculus.

%%%%%%%%%%%%%%%%%%%%%%%%%%%%%%%%%%%%%%%%%%%%%%%%%%%%%%%%%%%%%%%%%
\section{Relation to Bittner's operator calculus}
%%%%%%%%%%%%%%%%%%%%%%%%%%%%%%%%%%%%%%%%%%%%%%%%%%%%%%%%%%%%%%%%%
We can relate Bittner's operator calculus outlined in Appendix \ref{Appendix_BittnersOperatorCalculus} to the exterior derivative and homotopy operator. We expand the analogy that was introduced in \cite{KyciaPoincareCohomotopy}.

\begin{figure}
\centering
\xymatrix{   & 0       &   &  &   \\
L_{1}:=& \ar[u]_{\hat{d}} \mathcal{E} \ar@/^/[drr]^{H}  &  &   &  \\
L_{0} := & \ar[dr]_{\hat{d}} \mathcal{E} & \oplus & \mathcal{A} \ar@/^/[ull]^{d}  \ar[dl]^H   &  \\
     &  &   0   &
}
\caption{Operator calculus mapped to exterior calculus.}
\label{Fig.DecompositionOperatorCalculusExterior}
\end{figure}
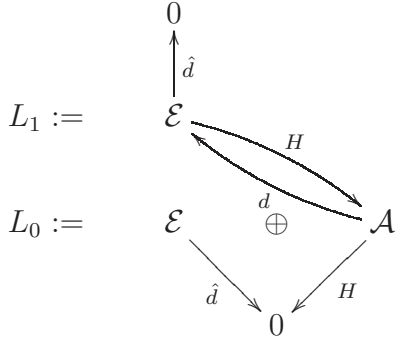

We define $L_{0}=\mathcal{E}\oplus\mathcal{A}$ and $L_{1}=\mathcal{E}$ as presented in Fig. \ref{Fig.DecompositionOperatorCalculusExterior}.
%%%%%%%%%%%
%%%%%%%%%%%
We have:
\begin{itemize}
 \item {$S:=d:L_{0}\rightarrow L_{1}$ - derivative with $ker(S)=ker(d)=\mathcal{E}$.}
 \item {$T:=H:L_{1}\rightarrow \mathcal{A}\subset L_{0}$ - integral.}
\end{itemize}
In addition, $ST|_{L_{1}}=dH|_{\mathcal{E}} = I$ since $dH$ is the projection operator onto $\mathcal{E}$.

In order to identify $s$ operator, we use homotopy invariance formula (\ref{Eq_homotopyFormula}) as
\begin{equation}
 Hd = I - \underbrace{(s_{x_{0}}^{*} + dH)}_{s},
\end{equation}
i.e.
\begin{equation}
 s:=\left\{
 \begin{array}{ccc}
  s_{x_{0}}^{*} & for & \Lambda^{0}(U) \\
  dH & for & \Lambda^{k}(U), \quad k>0.
 \end{array}
\right.
\end{equation}
The operator $s$ defined above is a projection ($s^2=s$) onto $ker(S)=ker(d)=\mathcal{E}$.

We will use the symbols $S$, $T$, and $s$ in these specific substitutions from exterior algebra.

We can now associate the operator $R=-A\wedge\_$ to the notion of abstract (non-commutative) logarithm in the sense that it has no zero divisors as it states the following Proposition.
%%%
\begin{Proposition}
 If $d\phi=R\phi$ then we have
 \begin{equation}
  (I-HR)\phi =0 \Rightarrow \phi=0.
 \end{equation}
\end{Proposition}
%%%
\begin{Proof}
If $\phi\in ker(R)$ then from $(I-HR)\phi=0$ we get that $\phi=0$.

If $\phi \not\in ker(R)$ then for the solution $\phi$ from (\ref{Eq.Solution_homogenous_k_gt_0}) we have $c = (I+H(A\wedge\_)\phi = (I-HR)\phi = 0$, and therefore again from (\ref{Eq.Solution_homogenous_k_gt_0}), $\phi=0$, as required.
\end{Proof}

The following Corollary instantiates a strict definition of formal notation of fraction in (\ref{Eq_inverison_of_1HA}).

\begin{Corollary}
 The operator $(I-HR)=(I+H(A\wedge\_))$ has no trivial zero divisors. Therefore, we can construct Mikusinki's ring of elements $\Lambda^{*}(U)$ and operators $\{ I, I-HR\}$, where $\phi \in \Lambda^{*}(U)$ is represented by $\frac{\phi}{I}$, and other elements are of the form $\frac{\phi}{I-HR}$.
\end{Corollary}

The following section provides a practical approach to solving the curvature equation.

%%%%%%%%%%%%%%%%%%%%%%%%%%%%%%%%%%%%%%
\section{Curvature equation}
%%%%%%%%%%%%%%%%%%%%%%%%%%%%%%%%%%%%%%
The curvature $F$ is a square of the differential operator $S-R$
\begin{equation}
 F:=(S-R)^{2}.
\end{equation}
We fix, as in the previous section, i.e., $S=d$ and $R=-A\wedge\_$, $s=dH+s_{x_{0}}^{*}$ so we obtain the curvature of $A$.

The curvature equation, due to its tensorial character, is an algebraic equation, however, we will focus on the solutions using the above machinery of (anti)exact forms to get a deeper insight into the structure of the solutions.

We can easily solve the curvature equation using the method of solving the covariant constancy equation. First, we solve the homogenous curvature equation in
%%%%%%%%%%%%%%%%%%%%%%%%%%%%%%%%%%
\begin{Theorem}
\label{Th_homogenous_curvature}
 In order to solve the homogenous curvature equation
 \begin{equation}
  F\phi = (S-R)^{2}\phi=0,
 \end{equation}
 with $\phi\in \Lambda^{k}(U,V)$ rewrite it as a coupled system
  \begin{equation}
  \left\{
  \begin{array}{cc}
   \phi_{2} := (S-R)\phi_{1}, & s\phi_{1}=0 \\
   (S-R)\phi_{2}=0, & s\phi_{2}=c_{2},
  \end{array}
  \right.
  \label{Eq_curvatureSystemHomogenous}
 \end{equation}
 for $\phi_{1}$ and $\phi_{2}$. To this solution one can add a solution to the first-order equation
 \begin{equation}
  (S-R)\phi=0, \quad s\phi=c_{1}\in ker(S).
 \end{equation}
From the first equation of (\ref{Eq_curvatureSystemHomogenous}) and Theorem \ref{Th_FullInhomogenous_parallelTransportEquation} one gets that the solution to (\ref{Eq_curvatureSystemHomogenous}) does not exist if $Hd\phi_{2} \not\in Im(A\wedge\_)$.
\end{Theorem}
%%%%%%%%%%%%%%%%%%%%%%%%%%%%%%%%%%

Next, the solution of the inhomogeneous curvature equation will be provided.
%%%%%%%%%%%%%%%%%%%%%%%%%%%%%%%%%%
\begin{Theorem}
 The solution of the inhomogeneous curvature equation
 \begin{equation}
  (S-R)^{2}\phi=J, \quad s\phi=c_{1}\in ker(S), \quad s(S-R)\phi=c_{2} \in ker(S),
 \end{equation}
 for $\phi \in \Lambda^{k}(U,V)$, $J \in \Lambda^{k+2}(U,V)$, $c_{1} \in \mathcal{E}^{k}(U,V)$, $c_{2} \in \mathcal{E}^{k+1}(U,V)$, $R\in \Lambda^{1}(U,End(V))$ is a linear combination of
\begin{itemize}
 \item {First order equation
 \begin{equation}
 \begin{array}{cc}
  (S-R)\phi=0 & s\phi=c_{1}.
 \end{array}
 \end{equation}
 The solution is provided by Theorem \ref{Th_homogenous_curvature}.
 }
 \item {Second order equation
 \begin{equation}
  (S-R)^{2}\phi=J, \quad s\phi=0, \quad s(S-R)\phi=c_{2}.
 \end{equation}
 It can be solved by replacing it with the first-order system of coupled equations:
  \begin{equation}
  \left\{
  \begin{array}{cc}
   \phi_{2} = (S-R)\phi_{1}, & s\phi_{1}=0 \\
   (S-R)\phi_{2}=J, & s\phi_{2}=c_{2},
  \end{array}
  \right.
  \label{Eq_curvatureInhomogenousSystem}
 \end{equation}
 with $\phi=\phi_{1}$ and $\phi_{2}$.

 From both equation (\ref{Eq_curvatureInhomogenousSystem}) applying successively Theorem \ref{Th_FullInhomogenous_parallelTransportEquation} one gets that the solution to (\ref{Eq_curvatureInhomogenousSystem}) does not exist if $Hd\phi_{2} \not\in Im(A\wedge\_)$ and if $HdJ \not\in Im(A\wedge\_)$.
 }
\end{itemize}

\end{Theorem}
%%%%%%%%%%%%%%%%%%%%%%%%%%%%%%%%%%
As pointed out, these solutions only have to exist if some constraints from the previous sections are met. However, the operator approach solves curvature equations as algorithmic as solving second order ODE.

%%%%%%%%%%%%%%%%%%%%%%%%%%%%%%%%%%%%%%
\section{Hodge duals}
\label{Section_Hodge_duals}
%%%%%%%%%%%%%%%%%%%%%%%%%%%%%%%%%%%%%%
When $U$ has a metric structure so the Hodge star $\star$ can be defined, we can use dual (anti)coexact decomposition for solving dual equations. The equation becomes
 \begin{equation}
  \delta\phi+A^{\sharp}\lrcorner\phi=0.
 \end{equation}
The term $A^{\sharp}\lrcorner\phi$ is understood as follows: $A^{\sharp}$, after choosing the base, is a matrix of vectors -  elements of $A$ matrix dualized by $\sharp$; then the operation $\lrcorner$ is both the matrix multiplication and insertion of elements (vectors) of $A^{\sharp}$ into elements of $\phi$.

By a similar argument, we must also put emphasis on elements of $ker(A^{\sharp}\lrcorner\_)$ in our considerations. The elements of $\mathcal{C}(U, V) \cap ker(A^{\sharp}\lrcorner\_)$ are \textit{coexact gauge modes} that introduce the nonuniqueness of the solution. Elements of the form $\mathcal{Y}(U,V) \cap ker(A^{\sharp}\lrcorner\_)$ are important in solving an inhomogenous equation and will be called \textit{anticoexact gauge modes}.

By dualizing Theorem \ref{Th_homogenous_solution} we have
%%%%%%%%%%%%%%%%%%%%%%%%%%%%%%%
\begin{Corollary}
 Solution of the equation
 \begin{equation}
  \delta\phi+A^{\sharp}\lrcorner\phi=0,
 \end{equation}
 where $\phi\in\Lambda^{k}(U,V)\setminus ker(A^{\sharp}\lrcorner\_)$, $A\in\Lambda^{1}(U,V)$ is given by
\begin{equation}
 \phi = \frac{1}{I+h(A^{\sharp}\lrcorner\_)}c = (I-h(A^{\sharp}\lrcorner\_)+h(A^{\sharp}\lrcorner(h(A^{\sharp}\lrcorner\_)))-\ldots)c,
\end{equation}
where $c\in ker(\delta)|_{U}\setminus ker(A^{\sharp}\lrcorner\_)$. The series is uniformly convergent to a continuous solution in the ball $B(x_{0},r)$ for the radius $r$ given by $r \frac{||A||_{\infty}}{k} =1,$, where the supremum norm is taken over the ball.

We can add to this solution coexact gauge mode $\mathcal{C}(U, V) \cap ker(A^{\sharp}\lrcorner\_)$. The space of solutions is, therefore, an affine space $\Lambda^{k}(U,V)\setminus ker(A^{\sharp}\lrcorner\_)$ over the vector space of coexact gauge modes.
\end{Corollary}
%%%%%%%%%%%%%%%%%%%%%%%%%%%%%%%

Likewise, from Theorem \ref{Th_nonhomogenous_solution_exactRHS} the inhomogeneous equation with coexacact RHS is provided by
\begin{Theorem}
  Solution of the equation
 \begin{equation}
  \delta\phi+A^{\sharp}\lrcorner\phi=J,
 \end{equation}
where $\phi\in\Lambda^{k}(U,V)\setminus ker(A^{\sharp}\lrcorner\_)$, $A\in\Lambda^{1}(U,V)$, $J\in\mathcal{C}^{k-1}(U,V)$ is given by
\begin{equation}
 \phi= \frac{1}{I+h(A^{\sharp}\lrcorner\_)}c +  \sum_{l=0}^{\infty} (-1)^{l} (h(A^{\sharp}\lrcorner\_))^{l} hJ
\end{equation}
where $c\in \mathcal{C}^{k}(U,V)\setminus ker(A^{\sharp}\lrcorner\_)$. The first term is a solution to the homogenous equation. The series is uniformly convergent to continuous solution in the ball $B(x_{0},r)$ for the radius $r$ given by $r \frac{||A||_{\infty}}{k} =1,$, where the supremum norm is taken over the closure of the ball.
\end{Theorem}

Finally, the solution of the covariant constancy equation with arbitrary inhomogeneity is given by
\begin{Theorem}
 Solution of the equation
 \begin{equation}
  \delta\phi+A^{\sharp}\lrcorner\phi=J,
 \end{equation}
where $\phi\in\Lambda^{k}(U,V)$, $A\in\Lambda^{1}(U,End(V))$, $J\in\Lambda^{k-1}(U,V)$ can be composed from three elements $\phi=\phi_{1}+\phi_{2}+\phi_{3}$, where $\phi_{1}$ is a solution of
\begin{equation}
 (\delta + A^{\sharp}\lrcorner )\phi_{1}=J_{c} -\delta(\phi_{2}+\phi_{3}),
\end{equation}
the $\phi_{2}$ is a solution of a constraint equation
\begin{equation}
 A^{\sharp}\lrcorner\phi_{2}=J_{y},
 \label{Eq_generalCosolutionConstraint}
\end{equation}
and $\phi_{3}$ is a solution of
\begin{equation}
 A^{\sharp}\lrcorner\phi_{3}=0,
\end{equation}
where $J_{c}=\delta hJ$ is the coecact part of $J$, and $J_{y}=h\delta J$ is the anticoexact part of $J$.

Note that from the constraint (\ref{Eq_generalCosolutionConstraint}) we get that the solution to the problem does not exist if $J_{y}\not\in Im(A^{\sharp}\lrcorner\_)$.
\end{Theorem}

We can also consider the square of the $\delta +A^{\sharp}\lrcorner\_$ operator, however, it can also be noticed that it is related to the results of the previous section, since by (\ref{Eq.HodgeDuality}) we have for $\alpha \in \Lambda^{k}(U,V)$ we have $(\delta + A^{\sharp}\lrcorner)\star \alpha=(-1)^{k+1}\star (d+\_\wedge A) \alpha$.

%%%%%%%%%%
\section{Cartan structure equations}
%%%%%%%%%
We now proceed to solve the Cartan structure equations \cite{Sternber}. The first equation
\begin{equation}
 d\theta + A\wedge \theta = \tau
\end{equation}
where $A\in \Lambda^{1}(U, End(V))$, and $\tau \in \Lambda^{deg(\theta)+1}(U, V)$ is related to a local form of a torsion tensor, is exactly inhomogeneous parallel transport equation considered before.

The second Cartan structure equation
\begin{equation}
dA+A\wedge A = F,
\label{Eq.SecondStructureEquation}
\end{equation}
where $F\in \Lambda^{2}(U, End(V))$ is the local form of the curvature, and the equation is for $A$. This equation will be the main part of this section. We have the following
%%%%%%%%%%%%%
\begin{Theorem}
For fixed $\alpha\in \mathcal{E}^{1}(U, End(V))$ the unique solution of (\ref{Eq.SecondStructureEquation}) is a limit of the following limit sequence
\begin{equation}
A_{0}=\alpha,~ A_{n} = \alpha -H(A_{n-1}\wedge A_{n-1}-F), ~ A = \lim_{n\rightarrow \infty}A_{n}.
\end{equation}
The limit exists for
\begin{equation}
 ||x-x_{0}||max(||\alpha\wedge \alpha ||_{\infty}, ||F||_{\infty}) < 2,
\end{equation}
and is continuous, i.e., $A \in \Lambda^{1,0}(U, End(V))$, with $dA$ being continuous form.
\end{Theorem}
%%%%%%%%%%%%%
The proof is rather standard usage of the theory of nonlinear integral equations to the exterior algebraic setup.
%%%%%%%%%%%%%
\begin{Proof}
 \textbf{Integral equation}: The first step is to bring (\ref{Eq.SecondStructureEquation}) into integral equation. Taking exterior derivative of both sides we get that $d(A\wedge A - F) =0$, i.e., $A\wedge A - F \in \mathcal{E}^{2}(U, End(V))$, and therefore, $A\wedge A-F = dH(A\wedge A - F)$. Substituting back to the original equation, we have
 \begin{equation}
  A = \alpha + H(A\wedge A-F),
  \label{Eq.SecondStructuralEq_IntegralForm}
 \end{equation}
for an arbitrary $\alpha \in \mathcal{E}^{1}(U,End(V))$. The equation (\ref{Eq.SecondStructuralEq_IntegralForm}) is a nonlinear integral equation for $A$.
\\
\textbf{Formal solution}: We solve (\ref{Eq.SecondStructuralEq_IntegralForm}) using iteration procedure:
\begin{equation}
 \begin{array}{c}
  A_{0}=\alpha, \\
  \vdots \\
  A_{n} = \alpha - H(A_{n-1}\wedge A_{n-1} -F), \\
  \vdots
 \end{array}
\end{equation}
Then, formally, the solution has the form
\begin{equation}
 A = \lim_{n\rightarrow\infty}A_{n}.
 \label{Eq.SecondCartanStructureEq_LimitSolution}
\end{equation}
\\
\textbf{Convergence and uniqueness}: In order to prove that the limit (\ref{Eq.SecondCartanStructureEq_LimitSolution}) converges to the continuous one-form, we reformulate it as a series:
\begin{equation}
 S = A_{0} + (A_{1}-A_{0})+(A_{2}-A_{1}) + \ldots.
\end{equation}
Using the norm (\ref{Eq.Norm_End_form}) we have
\begin{equation}
 \begin{array}{c}
  ||A_{0}||_{\infty} = ||\alpha||_{\infty}, \\
  ||A_{1} - A_{0}||_{\infty} = ||H(A_{0}\wedge A_{0} - F)||_{\infty} \leq \frac{1}{2}||x-x_{0}|| max(||A_{0}\wedge A_{0} ||_{\infty}, ||F||_{\infty}).
 \end{array}
\end{equation}
Moreover since the map $A\rightarrow A\wedge A$ has derivative $D (A\rightarrow A\wedge A)|_{A}(\gamma) = A\wedge \gamma + \gamma \wedge A$ for\footnote{We identify vectors space $\Lambda^{1}(U, End(V))$ and its tangent $T\Lambda^{1}(U, End(V))$.} $\gamma \in T\Lambda^{1}(U, End(V))$, so $||A_{n}\wedge A_{n} - A_{n-1}\wedge A_{n-1}||_{\infty} \leq 2||A_{n}-A_{n-1}||_{\infty}$. We therefore have
\begin{equation}
||A_{n}-A_{n-1}||_{\infty} \leq \frac{1}{2} ||x-x_{0}|| ||A_{n-1}-A_{n-2}||_{\infty} \leq \ldots \leq \left( \frac{||x-x_{0}||}{2}\right)^{n} max(||A_{0}\wedge A_{0} ||_{\infty}, ||F||_{\infty}).
\end{equation}
Then the series $||S||_{\infty}$ is uniformly convergent for
\begin{equation}
 \frac{1}{2}||x-x_{0}||max(||A_{0}\wedge A_{0} ||_{\infty}, ||F||_{\infty}) <1,
\end{equation}
from which results in the value of the radius of convergence. The series converges to a continuous function. Since $A$ also fulfills the original equation, so $dA$ is also continuous.

For the uniqueness of the solution, assume that there are two solutions $A$ and $\bar{A}$ of the integral equation. Then we have
\begin{equation}
max_{ij}||A_{ij}(x) - \bar{A}_{ij}(x)|| \leq ||H(A\wedge A - \bar{A}\wedge \bar{A})|| \leq ||x-x_{0}|| ||A-\bar{A}||_{\infty},
\end{equation}
which gives
\begin{equation}
 ||A-\bar{A}||_{\infty} \leq ||x-x_{0}|| ||A-\bar{A}||_{\infty},
\end{equation}
that it fulfilled for all $x\in U$ when $||A-\bar{A}||_{\infty} =0$.
\end{Proof}
%%%%%%%%%%%%%%%%%

We now pass to the simple scalar example.
\begin{Example}
Let us consider $\mathbb{R}^{2}$ with coordinates $x$, $y$. Let the curvature be $F = dx\wedge dy$.

For computing the connection, let us take as the starting value of the connection $A_{0}=dx$. Proceeding with the algorithm, we have,
\begin{equation}
 \begin{array}{c}
  A_{1} = dx + H(dx\wedge dy) = dx + \frac{1}{2}(xdy - ydx), \\
  A_{2} = dx + H(dx\wedge dy) = dx + \frac{1}{2}(xdy - ydx), \\
  \vdots
 \end{array}
\end{equation}
We observe that the derivations stabilize, and therefore,
\begin{equation}
 A = \lim_{n\rightarrow\infty} A_{n} = dx + \frac{1}{2}(xdy - ydx).
\end{equation}

One can note that the curvature $F$ can also arise from the other types of connections, e.g., $\bar{A}=xdy$ or $\hat{A}=-ydx$ plus an arbitrary exact part from $\mathcal{E}^{1}(U, End(V))$. Both of these connections have the same projector $Hd\bar{A}=Hd\hat{A}= \frac{1}{2}(xdy - ydx)$.
\end{Example}
%%%%%%%%%%%

In the next section we focus on the general nature of the exterior differential equations that contain $\dcov$ and its Hodge/metric dual.

%%%%%%%%%%
\section{General geometric-based differential equations}
%%%%%%%%%
The curvature is an example of a composition of $\dcov$ that is a tensor operator. We also have the Hodge dual to $\dcov$. In general, the composition of these operators is only a differential operator. We utilize the results of the preceding sections to solve the equation consisting of $\dcov$ for different connection one-forms and their Hodge duals, which will be called \textit{geometric-based differential operators}.
\begin{Definition}
For a one-forms $A_{i}\in \Lambda^{1}(U, End(V))$, $i\in\{1,\ldots, l_{1}\}$, and vector fields $X_{j}\in \Gamma(U)$, $\in \{1,\ldots l_{2}\}$ we define a covariant derivative operators associated with $A_{i}$
\begin{equation}
 D_{A_{i}}:=d+A_{i}\wedge\_:\Lambda^{k}(U,V)\rightarrow\Lambda^{k+1}(U,V),
\end{equation}
and dual covariant derivatives associated with vector fields $X_{i}$
\begin{equation}
\reflectbox{D}_{X_{i}}:=\delta+X_{i}\lrcorner\_:\Lambda^{k+1}(U,V)\rightarrow \Lambda^{k}(U,V).
\end{equation}
For a given set of indices, possibly some of them vanish, $P=(p_{1}, \ldots, p_{m})$ and $Q=(q_{1}, \ldots, q_{m})$ we define a geometric-based differential operator
\begin{equation}
D^{P,Q}:= D_{A_{p_{1}}}^{p_{1}}\reflectbox{D}_{X_{q_{1}}}^{q_{1}}\ldots D_{A_{p_{m}}}^{p_{m}}\reflectbox{D}_{X_{q_{m}}}^{q_{m}}:\Lambda^{k}(U,V)\rightarrow\Lambda^{k+\sum_{l=1}^{m}p_{l}-\sum_{l=1}^{m}q_{l}}(U,V).
\end{equation}
We set $D^{p}_{A}=\underbrace{D_{A}\circ \ldots\circ D_{A}}_{p}$, and similarly $\reflectbox{D}_{X}^{q}=\underbrace{\reflectbox{D}_{X}\circ \ldots \circ \reflectbox{D}_{X}}_{q}$. Moreover, $D_{A}^{0}=Id$ and likewise $\reflectbox{D}_{X}^{0}=Id$.

When $P=(1)$ and $Q=(0)$ (i.e., $D^{P,Q}=D_{A}$) or when $P=(0)$ and $Q=(1)$ (i.e., $D^{P,Q}=\reflectbox{D}_{X}$) we call $D^{P,Q}$ a first degree operator.

We can now define the (inhomogeneous) geometric-based differential equations
\begin{equation}
D^{P,Q}\phi = J,
\label{Eq_GeometricBasedDifferentialEquation}
\end{equation}
for $J\in \Lambda^{k+\sum_{l=1}^{m}p_{l}-\sum_{l=1}^{m}q_{l}}(U,V)$ with a solution $\phi \in \Lambda^{k}(U,V)$. When $J=0$ the equation is called homogenous.
\end{Definition}

As an example of geometric based differential operators we have $\delta d$ and $d\delta$, as well as, e.g., $d\delta + d(A\wedge\_)$ or $d\delta d$.

Utilizing Theorem \ref{Th_FullInhomogenous_parallelTransportEquation} and the idea of replacing the equation (\ref{Eq_GeometricBasedDifferentialEquation}) with the system of first degree equations we have
%%%%%%%%%%
\begin{Algorithm}
For solving (\ref{Eq_GeometricBasedDifferentialEquation}), replace it with the system
\begin{equation}
\left\{
\begin{array}{c}
 \reflectbox{D}_{X_{q_{m}}}\phi^{q_{m}}_{1}=\phi^{q_{m}}_{2} \\
 \reflectbox{D}_{X_{q_{m}}}\phi^{q_{m}}_{2}=\phi^{q_{m}}_{3} \\
 \ldots \\
 \reflectbox{D}_{X_{q_{m}}}\phi^{q_{m}}_{q_{m}-1}=\phi^{q_{m}}_{q_m} \\
 D_{A_{p_{m}}}\phi^{p_{m}}_{1}=\phi^{q_{m}}_{q_m} \\
 \ldots \\
 D_{A_{p_{m}}}\phi^{p_{m}}_{p_{m}-1}=\phi^{p_{m}}_{p_m}\\
 \ldots \\
 \ldots \\
 D_{A_{p_{1}}}\phi^{p_{1}}_{p_{1}}=\phi^{p_{1}}_{p_{1}-1}\\
 D_{A_{p_{1}}}\phi^{p_{1}}_{p_{1}} = J,
\end{array}
\right.
\end{equation}
for new forms $\phi^{p}_{q}$. Then solve the system starting from the last equation and proceed iteratively to the top using Theorem \ref{Th_FullInhomogenous_parallelTransportEquation}. The solution of the equation (\ref{Eq_GeometricBasedDifferentialEquation}) is $\phi= \phi^{q_{m}}_{1}$.

The initial data for all $\phi^{p}_{q}$ solutions are provided in terms of projection to the boundary conditions $s\phi^{p}_{q}=c^{p}_{q}$ for initial data $c^{p}_{q}$, and $s$ is either $s^{*}_{x_{0}}$ or $S_{x_{0}}=\star^{-1} s^{*}_{x_{0}}\star$ depending on if $\phi^{p}_{q}$ results from solving an equation containing $D_{A}$ or $\reflectbox{D}_{X}$ operators.

The possible obstruction during iterative solution results in the nonexistence of solution, as Theorem \ref{Th_FullInhomogenous_parallelTransportEquation} explains.
\end{Algorithm}
The above Algorithm provides a strict set of steps, similar to solving a higher-order ODE by reducing it to the system of first-order ODEs. The resemblance is apparent when one casts the results into Bittner's operator calculus approach.

\begin{Example}
As an example we consider the equation $\delta d\phi=0$ for $\phi\in\mathcal{E}(U,V)$. We rewrite it as a system of equations
\begin{equation}
\left\{
\begin{array}{c}
 d\phi_{1}=\phi_{2}\\
 \delta\phi_{2}=0,
\end{array}
\right.
\end{equation}
where $\phi=\phi_{1}$ is our solution. We first solve the second equation that gives $\phi_{2}=c \in \mathcal{C}(U,V)$ which is an arbitrary form. It is an exact gauge mode since there is no connection one-form, $A=0$. Then we substitute $\phi_{2}$ to the first equation obtaining the constraint $Hdc=0$ (i.e., $Hdc \in ker(A\wedge\_)=\{0\}$ by Theorem \ref{Th_FullInhomogenous_parallelTransportEquation}), and then the solution is $\phi=\phi_{1}=e + Hc$ for $e\in \mathcal{E}(U,V)$.

Likewise one can solve $d\delta \psi=0$ getting $\psi=f+hg$, where $g\in \mathcal{E}(U,V)$ with constraint $h\delta g=0$, and $f \in \mathcal{C}(U,V)$ is arbitrary.

The common solution of these two is the kernel of the Laplace-Beltrami operator $\Delta=d\delta + \delta d$. This imposes additional constraints on $\phi$ and $\psi$.
\end{Example}

%%%%%%%%%%
\section{Conclusions}
%%%%%%%%%
The formulas for inverting covariant exterior derivatives in a local star-shaped subset of a fibred set are provided. Using this prescription, the method of solving for the curvature equation is given. Moreover, the relation of the techniques developed here with Operational Calculus, especially with Bittner's calculus, is provided. Since Operational Calculus was invented to solve (linear) differential equations appearing in engineering as easily as algebraic equations, this link helps simplify notation and promotes an efficient way to make local calculations in differential geometry and their applications. Using (anti)(co)exact decomposition of forms allows solving equations involving covariant exterior derivatives as efficiently as standard ODEs.

\appendix

\section*{Acknowledgments}

This research was supported by the GACR grant GA22-00091S, the grant 8J20DE004 of the Ministry of Education, Youth and Sports of the CR, and Masaryk University grant MUNI/A/1092/2021. RK also thank the SyMat COST Action (CA18223) for partial support.

\section{Bittner's operator calculus}
\label{Appendix_BittnersOperatorCalculus}
Operational Calculus (e.g., \cite{Heaviside1, Heaviside2, Mikusinski, BergOperationalCalculus, DimovskiOperationalCalculus, Rolewicz}) was invented to deal with linear differential equations appearing in technical applications in an algebraic way. We are interested in the version of this calculus developed by R. Bittner \cite{Bittner1, Bittner2, BittnerBook1, BittnerBook1}. Since it is not well-recognized in differential geometry, we summarize the basic facts in this Appendix following \cite{Bittner1, Bittner2, BittnerBook1, BittnerBook1}.

We also acknowledge other similar approaches to the subject, especially to nonlinear problems, e.g., in \cite{Rosenkranz1}, and some use of skew-symmetry to operational calculus \cite{Rosenkranz2}, although they are using different axioms for derivative and integral. It can possibly be merged with the approach from this paper.

The derivative can be described in two ways: in a 'set-theoretic spirit' - describing how it acts on some elements of a set, or in a 'category-theory spirit' - describing how it interacts with other operators.
In calculus or differential geometry, a derivative (like covariant derivative) is defined by its linearity property and (possibly graded) Leibniz rule. On the contrary, the general idea for Bittner's operator calculus is to describe the derivative using its interplay with an integral operator. It generalizes the following rules of the elementary calculus
\begin{equation}
 \frac{d}{dx}\int_{x_{0}}^{x}f(t)dt = f(x), \quad \int_{x_{0}}^{x}\frac{df}{dt}(t)dt = f(x)-f(x_{0}),
\end{equation}
for $f\in C^{\infty}(\mathbb{R})$.

\begin{Definition} \cite{BittnerBook1, BittnerBook2, Bittner1, Bittner2}
 For linear spaces $L_{0}$ and $L_{1}$ we define linear operators
 \begin{itemize}
  \item {$S:L_{0}\rightarrow L_{1}$ - abstract derivative;}
  \item {$T_q:L_{1}\rightarrow L_{0}$ - abstract integral parametrized by $q\in ker(S)\subset L_{0}$;}
  \item {$s:L_{0}\rightarrow ker(S)\subset L_{0}$ - projection/limit condition;}
 \end{itemize}
 that fulfills
 \begin{equation}
  ST= I, \quad TS = I - s.
 \end{equation}

 Elements of $ker(S)$ are called constants (of $S$).
\end{Definition}
The $T$ is the right inverse of $S$, however, if $s\neq 0$, there is a defect \cite{Rolewicz} preventing $T$ to being also the left inverse of $S$. The operator $s$ is a projection operator, so $s^{2}=s$. We can extend the calculus to higher derivatives by the chain of linear spaces $\{L_{i}\}_{i=0}^{\infty}$ and derivatives/integrals acting between them: $L_{0}\xrightleftharpoons[T]{S} L_{1}\xrightleftharpoons[T]{S} L_{2}\xrightleftharpoons[T]{S} \ldots$, however, we do not need it since we will use this calculus here to the exterior derivative that is nilpotent, and so this chain will terminate after the second term.

This operator calculus can be defined even more abstractly on groups by
\begin{Definition} \cite{Bittner2, BittnerBook2}
 We define for groups $G_{0}$, $G_{1}$
 \begin{itemize}
  \item {$S:G_{1}\rightarrow G_{0}$, $S$ being group homomorphism - abstract derivative;}
  \item {$T:G_{0}\rightarrow G_{1}$ - abstract integral for $S$;}
  \item {$s:G_{1}\rightarrow ker(S)\subset G_{1}$ - projection/limit condition;}
 \end{itemize}
 with properties
 \begin{equation}
  \forall_{y\in Y}~ STy=y, \quad \forall_{x\in X}~ x=s(x) \circ TS(X).
 \end{equation}
\end{Definition}
As an example consider $G_{0}=(\mathbb{R}_{+}\sum\{0\},1,\cdot)$, $G_{1}=(\mathbb{R},1,\cdot)$ with $Sx=x^{2}$, $Tx=\sqrt{x}$, $s(x)=sgn(x)$, with obvious conditions $STx=x$, $y=s(y)TSy$.

For vector spaces treated as abelian groups $G_{i}=(V,0,+)$, $i=1,2$ we restore the previous definition. We will not need this abstract definition in what follows, however, it can be useful in the case of group manifold.

The abstract logarithm map is defined, as usual \cite{NapierLogarithm}, as a solution to some (abstract) differential equation. To this end, consider an endomorphism $R:L_{i}\rightarrow L_{i}$ for $i=0,1$ that commutes with $S$ and $s$. One can prove \cite{BittnerBook1, BittnerBook2} that it also commutes with $T$. We can define
\begin{Definition} \cite{BittnerBook1, Bittner2}
 The endomorphism $R$ that commutes with $S$ and $s$ is called \textit{logarithm} if fulfills
 \begin{equation}
  (I-TR)f=0 \Rightarrow f=0
  \label{Eq_logarithm_definition}
 \end{equation}
for $f\in L_{0}$. It is equivalent to demanding that
\begin{equation}
 \left\{ \begin{array}{c}
          Sf=Rf \\
          sf=0
         \end{array}
\right. \Rightarrow f=0.
\end{equation}
\end{Definition}
We will use the condition (\ref{Eq_logarithm_definition}) to generalize logarithm to the case when $R$ does not commute with $S$ and $s$, i.e., for operator calculus on exterior algebra.

The final definition from Operational Calculus we need is the definition of Mikusinki's ring (commutative ring without zero divisors). It was originally formulated for the convolutional ring of functions \cite{Mikusinski} as an alternative approach to distribution theory (it states a strict base of formal Heaviside's calculus \cite{Heaviside1, Heaviside2}). We can define it as follows
\begin{Definition}\cite{BittnerBook1, BittnerBook2, Mikusinski}
 Let $\pi(X)$ be the commutative group of endomorphisms of $X$ defined by
 \begin{itemize}
  \item {$(U_{1}U_{2})U_{3}=U_{1}(U_{2}U_{3})$}
  \item {$U_{1}U_{2}=U_{2}U_{1}$}
  \item {No zero divisors: $Uf=0 \Rightarrow f=0$}
 \end{itemize}
 for $U_{1},U_{2},U_{3}, U \in \pi(X)$ and $f\in X$.

 Then we can define fractions
 \begin{equation}
  \xi = \frac{f}{U}
 \end{equation}
as an equivalence relation: $\frac{f}{U}=\frac{g}{V} \Leftrightarrow Vf=Ug$. This defines the ring with the obvious addition of fractions and multiplication of a fraction by numbers and operators.
\end{Definition}
In the paper we define Mikusinski's ring for the group of operators containing one operator only.

To make a connection with differential geometry, we define a curvature $F$ as a square of some operator that defines a logarithm: $F=(S-R)^{2}$. The additional requirement expected in differential geometric applications is that the curvature should be an algebraic operator, i.e., not contain derivative acting on the argument. Therefore we have an obvious
\begin{Proposition}
 For an operator that defines logarithm, i.e., $S-R$ defines an algebraic curvature, when $S$ fulfills graded Leibniz rule, is nilpotent $S^{2}=0$, and $R$ has an odd grade.
\end{Proposition}
\begin{Proof}
 The proof is a straightforward computation
 \begin{equation}
  (S-R)^{2}\phi = -S(R\phi)-RS\phi+R^{2}\phi = -(S(R)-R^{2})\phi = F\phi,
 \end{equation}
where $F$ is a curvature operator that does not act with $S$ on $\phi$, i.e., it is an algebraic operator, as required.
\end{Proof}
%%%%%%%%%%%%%%%%%%%%%%%%%%%%%%%%%%%

%%%%%%%%%%%%%%%%%%%%%%%%%%%%%%%%%%%%%
\section{Bittner's operator calculus as category}
\label{Appendix_BittnersOperatorCalculusAsCategory}
%%%%%%%%%%%%%%%%%%%%%%%%%%%%%%%%%%%%%
In this appendix the category of Bittner's operator calculus is formulated, elevating it to the modern mathematical tool.

The category consists objects $B(L_{0},L_{1},S, T, s)$ that contains the following data
\begin{itemize}
 \item {$L_{0},L_{1}$ - linear spaces}
 \item {$S:L_{0}\rightarrow L_{1}$ - abstract derivative with $ker(S)\subset L_{0}$}
 \item {$T_q:L_{1}\rightarrow L_{0}$ - abstract integral indexed by $q\in ker(S)$}
 \item {$s_{q}:L_{0}\rightarrow ker(S)$ - projection}
\end{itemize}
with properties
\begin{equation}
 ST_{q}=I, \quad T_{q}S=I-s_{q}.
\end{equation}

The morphism between two such objects $B(L_{0},L_{1},S, T, s)$ and $\bar{B}(\bar{L_{0}},\bar{L_{1}},\bar{S}, \bar{T}, \bar{s})$ consists of two isomorphisms of linear spaces
\begin{equation}
 \phi:L_{1}\rightarrow \bar{L_{1}}, \quad \psi: L_{0}\rightarrow \bar{L_{0}},
\end{equation}
with the transformations of operators \cite{Bittner1, Bittner2, BittnerBook1, BittnerBook2}
\begin{equation}
 \bar{S}=\psi S\phi^{-1}, \quad \bar{T}=\phi T \psi^{-1}, \quad \bar{s}=\phi s \phi^{-1}.
\end{equation}

%%%%%%%%%%%%%%%%%%%%%%%%%%%%%%%%%%%%%
\section{Operator-valued connection}
%%%%%%%%%%%%%%%%%%%%%%%%%%%%%%%%%%%%%
In this part we discuss the application of the inversion of covariant derivative for operator-valued covariant derivative (see \cite{EdelenOperatorValuedConnection} or chapter 10 of \cite{EdelenExteriorCalculus}). The space is constructed from the product of Euclidean $E_{n}$ and the vector space $\mathbb{R}^{N}$:
\begin{equation}
 K:=E_{n}\times \mathbb{R}^{N},
\end{equation}
onto which acts an $r$-dimensional Lie group $G$ and its Lie algebra $\mathfrak{g}$. Enlarging the space to
\begin{equation}
 \mathcal{G}:= G\times K,
\end{equation}
we consider $\Lambda(\mathcal{G})$. Let Lie algebra is spanned by $\{V_{a}\}_{a=1}^{r}$, then the action of $V_{a} \in \mathfrak{g}$ on $\omega\in\Lambda(\mathcal{G})$ is given by the Lie derivative
\begin{equation}
 (V_{a}, \omega) \rightarrow \mathcal{L}_{V_{a}}\omega.
\end{equation}
It can be integrated into the action of the group by
\begin{equation}
 \omega\rightarrow \exp(u^{a}\mathcal{L}_{V_{a}})\omega,
\end{equation}
where $u^{a}$ are the coordinates on $G$. Then we can define operator-valued covariant derivative on $\mathcal{G}$ by
\begin{equation}
 D:= d + \Gamma , \quad \Gamma:=W^{a}\wedge\mathcal{L}_{V_{a}}, \quad W^{a} \in \Lambda^{1}(\mathcal{G}).
\end{equation}

Since the exterior derivative $d$ commutes with the Lie derivative $\mathcal{L}$, and the Lie derivative does not alter the grading of differential form, so all the results presented above are valid by replacing $A\wedge\_$ by $\Gamma$.

%%%%%%%%%%%%%%%%%%%%%%%%%%%%%%%%%%%%%
\section{Integral equations}
\label{Appendix_Integral_equations}
%%%%%%%%%%%%%%%%%%%%%%%%%%%%%%%%%%%%%

In this section we provide an alternative formulation in terms of integral equations generated by the homotopy operator $H$.

Let us consider the equation
\begin{equation}
 d\phi+A\wedge\phi = J_{e}, \quad A\in\Lambda^{1}(U), \quad J_{e}\in \mathcal{E}.
\end{equation}
Applying $d$ we get that $A\wedge \phi \in \mathcal{E}(U,V)$. Therefore, since $dH(A\wedge\phi)=A\wedge \phi$ and $dHJ_{e}=J_{e}$, so
\begin{equation}
 d(\phi+H(A\wedge\phi)-J_{e})=0,
\end{equation}
and therefore, we must solve the integral equation
\begin{equation}
 \phi+H(A\wedge\phi)=J_{e}+d\alpha,
\end{equation}
for an arbitrary $\alpha\in\Lambda(U,V)$. Now, introducing a real nonzero parameter $\lambda$, we get
\begin{equation}
 \phi+\lambda H(A\wedge\phi)=J_{e}+d\alpha.
 \label{Eq_IntegralEquationJExact}
\end{equation}
Substituting a Neumann series \cite{IntegralEquations}, $\phi=\sum_{l=0}^{\infty} \lambda^{l}\phi_{l}$, we restore the result of Theorem \ref{Th_nonhomogenous_solution_exactRHS}.

We also want to acknowledge that the methods presented here can be easily extended to solve the Riemann-Graves integral equation
\begin{equation}
 \phi = I +H(\phi \Gamma),
\end{equation}
for an unknown $\phi$ being a matrix-valued form and $\Gamma$ a matrix-valued $1$-form. The equation (\ref{Eq_IntegralEquationJExact}) can be seen as a generalized Riemann-Graves equation. It has applications in gauge theory \cite{RiemannGravesIntegralEquation, EdelenExteriorCalculus}.

%%%%%%%%%%%%%%%%%%%%%%%%%%%%%%%%%%%%%
\section{Antiexact inhomogeneity and introduction of connection form}
\label{Appendix_Antiexact_inhomogenity}
%%%%%%%%%%%%%%%%%%%%%%%%%%%%%%%%%%%%%
We now want to focus on the relation between the antiexact type of inhomogeneity in (\ref{Eq_nonhomogenous_covariant_equation}) and the term with the connection form $A$. Since this Appendix serves as a motivation and informal discussion, we focus on scalar-valued differential forms on $U$, however, the ideas can be easily extended to fibered spaces that ale local models of vector and associated vector bundles.

Start with a simpler equation of type (\ref{Eq_nonhomogenous_covariant_equation}) with $A=0$, i.e.,
\begin{equation}
d\phi = J,
\label{Eq_dphi=J}
\end{equation}
where $\phi, J \in \Lambda(U)$.
The standard argument for solving it \cite{EdelenExteriorCalculus} is as follows: by taking the exterior derivative of both sides, one gets consistency condition $dJ=0$, i.e., $J=dHJ \in \mathcal{E}(U)$. Then we have $d(\phi-HJ)=0$, that is $\phi = c + HJ$ for some arbitrary $c\in \mathcal{E}(U) = ker(d)$. This is precisely the same form as for a solution of a first-order ODE\footnote{For an ODE: $\frac{dx(t)}{dt}=f(t)$ for some $f \in C^{\infty}(\mathbb{R})$, namely, $x(t) = C + \int f(t) dt$, where $C\in \mathbb{R}=ker\left(\frac{d}{dt}\right)$.}. Therefore, $c$ can be considered as a constant of integration for an 'integral' $H$, or some kind of an 'initial condition', which is precisely the setup that fits into the definition of an abstract ODE within Bittner's operator calculus.

Now let us consider the case when there is a nonzero antiexact part in $J$, i.e., when $J=J_{e}+J_{a}$ for $J_{e}\in \mathcal{E}(U)$, $J_{a}\in \mathcal{A}(U)$ and $J_{a}\neq 0$. Then taking exterior derivative of both sides of (\ref{Eq_nonhomogenous_covariant_equation}) we get that $dJ_{a}=0$, i.e., $J_{a}\in \mathcal{E}(U)\cap \mathcal{A}(U)=\emptyset$, a contradiction. Therefore to balance the antiexact inhomogeneity of $J_{a}$ we must provide some additional term, e.g., $A\wedge \phi$ to obtain equation \ref{Eq_nonhomogenous_covariant_equation}), whose solutions were analyzed in subsection \ref{Subsection_Inhomogenous_equation}. One can therefore interpret antiexact inhomogeneity, $J_{a}$, that forces us to introduce the additional term with an 'external quantity' -- one-form $A$ -- as a 'non-autonomous term'\footnote{In analogy to the ODE $\frac{dx}{dt}=f(x,t)$, where the $t$ dependence of $f$ indicates that the system described by the ODE is part of a larger system and 'interacts' with it. We can also plot an analogy to Classical Field Theory \cite{Thirring, ConnectionsInPhysics}, where a system described by a field $\phi$ interacts with an external system by a 'current' $J$, which induces the new field $A$ that mediates in this interaction \cite{BennTucker, GagueTheories}.}.

One can also note that there is another way of augmenting (\ref{Eq_dphi=J}) for $J$ having a nonzero antiexact part. For instance, we can construct the following equation:
\begin{equation}
d\phi + H\beta = J_{e}+J_{a},
\end{equation}
for unknown $\phi$ and $\beta$. Then the system decouples into two independent equations
\begin{equation}
d\phi = J_{e}, \quad H\beta = J_{a},
\end{equation}
with solutions
\begin{equation}
\phi = c + HJ_{e}, \quad \beta = a + dJ_{a},
\end{equation}
for arbitrary $c\in\mathcal{E}(U)$ and $a\in \mathcal{A}(U)$. We must provide additional relation between $\phi$ and $\beta$ for coupling these equations. For obtaining (\ref{Eq.FullInhomogenous_CovariancyConstantEquation}), we can impose $H\beta =A\wedge \phi$, that is an algebraic equation for $A$, namely,
\begin{equation}
HdJ_{a}=A\wedge (c + HJ_{e}).
\end{equation}
That shows that the inhomogeneity $J$ provides $A$. By taking the exterior derivative and using $dHd=d$, one gets
\begin{equation}
 dJ_{a} = F\wedge (c + HJ_{e}) - A\wedge dHJ_{e} - A\wedge A \wedge (c+HJ_{e}),
\end{equation}
where $F=dA+A\wedge A$ is the curvature of $A$. This approach can be treated as an inverse problem -- fixing $J$, find $A$ such that $\dcov\phi=J$ is consistent.

Similar ideas can be applied to Hodge-dualized concepts on fibered sets with the Riemann metric by replacing $d \rightarrow \delta$ and $H \rightarrow h$. The equation
\begin{equation}
 \delta \phi = J
 \label{Eq_deltaPhi=J}
\end{equation}
is solved by
\begin{equation}
 \phi = c + hJ,
\end{equation}
where $c \in \mathcal{C}(U)$ and the consistency condition $\delta J=0$ states that $J\in \mathcal{Y}(U)$. When $J=J_{co}+J_{ac}$, for $J_{co}\in\mathcal{C}(U)$ and $J_{ac}\in\mathcal{Y}(U)$, $J_{ac}\neq 0$, we get contradiction unless we amend (\ref{Eq_deltaPhi=J}) with additional term. One of the possibilities is adding $A^{\sharp}\lrcorner \phi$, which leads to the equation analyzed in Section \ref{Section_Hodge_duals}.

%%%%%%%%%%%%%%%%%%%%%
\section{Regularity considerations}
\label{Section_RegularityConsiderations}
%%%%%%%%%%%%%%%%%%%%%
The obtained solution is a continuous form with continuous exterior derivative. It is a result of the fact that the uniform convergence of series of smooth functions converges to a continuous function.

When $A=0$ in the equation $\dcov \phi =0$ then the resulting equation $d\phi=0$ can have smooth solution, and therefore in this limit case we expect that the both series and its higher derivatives converges uniformly for arbitrary $x \in U$. We therefore expect that the radius of convergence of derivatives of (\ref{Eq.Solution_homogenous_k_gt_0}) is controlled by the 'largeness' of $A$.

Introduce for $\phi \in \Lambda^{*,l}(U,V)$ the norm
\begin{equation}
 || \phi ||_{l,\infty} = max_{|\alpha| \leq l}||D^{\alpha}\phi||_{\infty},
\end{equation}
where $D^{\alpha} = \frac{\partial^{|\alpha|}}{\partial x_{1}^{\alpha_{1}}\ldots \partial x_{n}^{\alpha_{n}}}$, for a multiidex $\alpha=(\alpha_{1},\ldots, \alpha_{n})$, and $|\alpha|=\alpha_{1}+\ldots+\alpha_{n}$. The supremum $||\_||_{\infty}$ is taken on some fixed compact set.

For constructing the solution (\ref{Eq.Solution_homogenous_k_gt_0}) with the coefficients from the class $C^{l}$, one have to provide uniform convergence of the solution $\phi$ and its derivatives up to $l$-th order, for which the norm $||\_||_{l,\infty}$ is the most suitable. We have for $|\alpha|\leq l$ that
\begin{equation}
\begin{array}{c}
     ||D^{\alpha}H(A\wedge\omega) ||_{l,\infty} \leq \frac{l^{2}}{l+k-1}  ||A||_{l,\infty} ||\omega||_{l,\infty} + \frac{l+1}{l+k}||x-x_{0}|| ||A||_{l,\infty}||\omega||_{l,\infty} \leq \\
     \leq \frac{1}{l+k-1}||A||_{l,\infty}||\omega||_{l,\infty}(l^{2} + (l+1)||x-x_{0}||)
\end{array}
\end{equation}
Therefore, the derivatives of order $l$ of terms of (\ref{Eq.Solution_homogenous_k_gt_0}) form uniformly convergent power series when
\begin{equation}
\frac{1}{l+k-1}||A||_{l,\infty}(l^{2} + (l+1)||x-x_{0}||) < 1.
\end{equation}
The inequality determines the boundary of uniform convergence. The prominent consequence is the fact that when $||A||_{l,\infty}$ is small then $||x-x_{0}||$ can be large, as intuition suggests.

The essential problem in obtaining higher regularity of the solution is a lack of theorem that relates uniformly convergent series of smooth functions to the smooth function. However, on the complex domain, uniformly convergent series of analytic functions converges to the analytic function. Therefore, when one focus on an a complex manifold $M$ with a star-shaped subset $U$, where all differential forms have analytic coefficients, then the solution (\ref{Eq.Solution_homogenous_k_gt_0}) can uniformly convergent to an analytic differential form.

%%%%%%%%%%%%%%%%%%%%%%%%%%%%%%%%%%%%%%%%%%%%%%%%%%%%%%%%%%%%%%%%%%%%%%%%%%%
%bibliografia
%\newpage
%%%%%%%%%%%%%%%%%%%%%%%%%%%%%%%%%%%%%%%%%%%%%%%%%%%%%%%%%%%%%%%%%%%%%%%%%%%

%%%%%%%%%%%%%%%%%%%%%%%%%%%%%%%%%%%%%%%%%%%%%%%%%%%%%%%%%%%%%%%%%%%%%%%%%
%koniec bibliografii
%%%%%%%%%%%%%%%%%%%%%%%%%%%%%%%%%%%%%%%%%%%%%%%%%%%%%%%%%%%%%%%%%%%%%%%%%%
%%%%%%%%%%%%%%%%%%%%%%%%%%%%%%%%%%%%%%%%%%%%%%%%%%%%%%%%%%%%%%%%%%%%%%%%%%%

%KONIEC
%%%%%%%%%%%%%%%%%%%%%%%%%%%%%%%%%%%%%%%%%%%%%%%%%%%%%%%%%%%%%%%%%%%%%%%%%

\end{document}